\documentclass[12pt]{article}
\usepackage{float} 
\usepackage{amsmath}
\usepackage{amssymb}
\usepackage{a4}
\usepackage{theorem}
\usepackage[mathscr]{eucal}
\usepackage{multirow}
\newcommand{\cn}{\mathop{\rm :}}
\newcommand{\CC}{\mathbb{C}}
\newcommand{\NN}{\mathbb{N}}
\newcommand{\PP}{\mathbb{P}}
\newcommand{\QQ}{\mathbb{Q}}
\newcommand{\RR}{\mathbb{R}}
\newcommand{\ZZ}{\mathbb{Z}}
\newcommand{\tC}{\tilde{C}}
\newcommand{\tE}{\tilde{E}}

\newcommand{\tX}{\tilde{X}}

\newcommand{\sL}{\mathscr{L}}
\newcommand{\sS}{\mathscr{S}}
\newcommand{\cO}{{\cal O}}
\newcommand{\oX}{\overline{X}}
\newcommand{\ol}{\overline}
\newcommand{\la}{\lambda}
\def\lh({\/{\rm (}}
\def\rh){\/{\rm )}}
\newcommand{\kod}{\kappa}
\newcommand{\mult}{\operatorname{mult}}
\newcommand{\piast}{\pi^\ast}

\newcommand{\smtX}{\smash{\tX}}
\newcommand{\brtX}{\left(\smtX\right)}
\newcommand{\KtX}{K_{\smtX}}
\newcommand{\KoX}{K_{\oX}}
\newcommand{\OoX}{\cO_{\oX}}
\newcommand{\KXs}{K_{X_s}}
\newcommand{\coneoX}{c_1\left(\smash{\oX}\right)^2}
\newcommand{\conetX}{c_1\brtX^2}
\newcommand{\ctwooX}{c_2\left(\smash{\oX}\right)}
\newcommand{\ctwotX}{c_2\brtX}
\newcommand{\etX}{e\brtX}
\newcommand{\eXs}{e\left(X_s\right)}
\newcommand{\pgtX}{p_g\left(\smash{\tX}\right)}
\newcommand{\pgXs}{p_g\left(X_s\right)}
\newcommand{\ptwotX}{P_2\left(\smash{\tX}\right)}
\newcommand{\qtX}{q\left(\smash{\tX}\right)}
\newcommand{\qXs}{q\left(X_s\right)}
\newcommand{\chitX}{\chi\left(\cO_{\smash{\tX}}\right)}
\newcommand{\chiXs}{\chi\left(\cO_{X_s}\right)}
\newcommand{\bzerotX}{b_0\left(\smash{\tX}\right)}
\newcommand{\bzeroXs}{b_0\left(X_s\right)}
\newcommand{\bonetX}{b_1\left(\smash{\tX}\right)}
\newcommand{\btwotX}{b_2\left(\smash{\tX}\right)}
\newcommand{\btwoXs}{b_2\left(X_s\right)}
\newcommand{\bthreetX}{b_3\left(\smash{\tX}\right)}
\newcommand{\honeonetX}{h^{1,1}\left(\smash{\tX}\right)}
\newcommand{\honeoneXs}{h^{1,1}\left(X_s\right)}
\newcommand{\mtd}{\mu_3\left(d\right)}
\newcommand{\KtC}{K_{\smash{\tC}}}
\newcommand{\kodtX}{\kod{(\tX)}}
{\theorembodyfont{\rm}
    \newtheorem{definition}{Definition}[section]
    \newtheorem{remark}[definition]{Remark}
    
}
\newtheorem{lemma}[definition]{Lemma}
\newtheorem{corollary}[definition]{Corollary}
\newtheorem{proposition}[definition]{Proposition}
\newtheorem{theorem}[definition]{Theorem}
\newcommand{\proof}{{\noindent\it Proof:\enspace }}
\newcommand{\proofend}{$\quad\square$\medskip\par}
\newcommand{\bibauthor}[2]{{\it{#2}~{#1},}}
\newcommand{\bibtitlea}[1]{#1,}            
\newcommand{\bibtitleb}[1]{{\it #1},}      
\newcommand{\bibcompany}[1]{#1}

\newcommand{\bibyear}[1]{(#1),}
\newcommand{\bibyearx}[1]{(#1)}
\newcommand{\bibjournal}[2]{#1\ {\bf #2}}
\newcommand{\bibpages}[2]{{#1}--{#2}}
\newcommand{\bibend}{.}

\begin{document}

\title{\bf Surfaces with triple points}
\author{Stephan Endra\ss\ \and Ulf Persson \and Jan Stevens
       \thanks{Partially supported by the 
               Swedish Natural Science Research Council (NFR)}}
\maketitle

\begin{abstract}
\noindent
In this paper we compute upper bounds for the 
number of ordinary triple points
on a hypersurface in $\PP^3$ and give a complete classification for 
degree six
(degree four or less is trivial, and five is elementary). But the real 
purpose is to 
point out the intricate geometry of examples with many triple points, and
how it fits with the general classification of surfaces.
\end{abstract}

\section*{Introduction}
The problem of finding 
the maximal number of simple singularities on projective hypersurfaces
has attracted a lot of attention in the last twenty years.
In this paper we study surfaces with a simple kind of
non-simple singularities, namely ordinary triple points.
In contrast to the case of simple surfaces singularities 
(classified by DuVal as those which `do not affect the 
conditions of adjunction' \cite{duval})
the invariants of the surface and its type in the classification
of surfaces may change.
Normal surfaces with higher singularities provide 
interesting examples of surfaces found, so to speak,
in our back-yard.

We warm up by looking at quintics with isolated triple points.
Their analysis is very elementary, 
but yields examples which nicely illustrate many aspects of the general theory
of surfaces. Quintics with many triple points
were to our knowledge first studied by Gallarati \cite{gallarati}.

For higher degree surfaces the location of the triple points
will matter very much, and the problem of finding the maximal number
becomes very difficult. We derive several
bounds on the number of triple points. We found
an example of a septic with a high number (16) of triple points,
which is one short of our upper bound. For sextics we give a
classification, which takes up the main part of this paper.

Our research started out as search for sextics with many triple points.
In particular, a sextic with 11 triple points would be very interesting. 
In fact, given the right configuration lying on a quadric, it would furnish 
a birational Abelian surface. However, we could not find such a surface; 
the putative example turned out to be a triple cover of a quadric. Then it 
was easily remarked that
11 triple points are {\it a priori} impossible.
Instead we found many examples with 9 triple points
and general arguments allowing to rule out possibilities.
The successful construction of an example raises the question how special
the construction is: can it be generalised? This can be decided
by infinitesimal methods. Let $\Sigma^d_\nu$ be the stratum of 
surfaces with $\nu$ ordinary triple points in the
parameter space of all surfaces of degree $d$. 
A lower bound for  the dimension of $\Sigma^d_\nu$ in the point
representing an explicit example $X$ is 
the number of moduli in the construction
plus $15$ from coordinate
transformations --- the stabiliser of the point configurations being discrete
for large $\nu$. 
An upper bound is given by the dimension of the Zariski tangent space.
If these dimensions are equal we know that the stratum is smooth of the given
dimension in that point, and the construction gives the general
element of the component of $\Sigma^d_\nu$.
To compute the Zariski tangent space we 
have to determine which polynomials of degree $d$ induce
equisingular deformations of the singular points, which
for each specific example can easily be computed with
{\it Macaulay} \cite{bayer}. 

The clue to classifying sextics with many triple points is
the study of exceptional curves of the first kind on the minimal
resolution. It turns out that only a few different cases can occur.
For nine triple points, we find three families
of $K3$ surfaces
(theorem \ref{theorem:K3}) and two families of properly
elliptic surfaces (theorem \ref{theorem:elliptic}).
In the $K3$ cases we find for each possible configuration
of the nine points a pencil of sextics of the form 
$\alpha g+\beta q^3$, where $q$ defines the (unique) quadric through
the nine points, and $g$ is a reducible sextic. In the other two cases
we find even a net of sextics, again containing $q^3$ and reducible sextics.
Constructing reducible sextics with triple points is not so
difficult.

Regarding the existence of 
a sextic $\{g=0\}$ with 10 triple points we first observe that the
pencil $\alpha g+\beta q^3$, where $\{q=0\}$ passes through nine of them,
falls into one of our five families of sextics with nine triple points.
Assuming that an element of such a family has a tenth triple point
gives conditions on the coefficients. The resulting equations
are much to difficult to solve. By imposing extra symmetry
we have been able to reduce the number of variables
and equations, while keeping at least a one-parameter family of solutions.

This paper is organised as follows.
First we study quintics with triple points.
The next section describes the birational invariants of our surfaces.
In the third section several bounds for the number of triple points
on a surface of degree $d$ are given. Then we identify
the exceptional curves.
The next section is the main part of this paper.
First we study the exceptional curves on a
sextic with triple points and find only a few
different cases, depending on the geometric genus
(corollary \ref{proposition:curves}).
Based on the geometry of the exceptional curves,
we give an overview over sextic surfaces with
$0$, \ldots, $10$ ordinary triple points.
For convenience, a
list of sextics with triple points and their invariants
is given (theorem \ref{theorem:sextics}).
In the last section a septic
with 16 ordinary triple points is constructed.
\section{Quintics}
A smooth quintic in $\PP^3$ is one of the simplest examples of
a surface of general type. Its Chern invariants are
given by $c_1^2=5$, $c_2=55$ and thus $\chi=4$. Furthermore it has the 
Hodge invariants
$p_g=4$, $q=0$ and $h^{1,1}=45$
(we give general formulas in the next section). 
The adjoint system consists of planes,
thus is $4$-dimensional (illustrating $p_g=4$). 

If nodes are allowed, or more generally rational double points,
nothing happens to the invariants, as the corresponding resolutions
are smooth deformations of smooth quintics. It is however an
interesting question to try and decide what 'bouquets' of rational
double points can be imposed. The maximal number of nodes is
31 (see below), 
and one may impose five $E_8$ singularities ($w^5=F(x,y,z)$ where $F$ is
a plane quintic with five cusps) which is somewhat short of the
maximal number  $44=h^{1,1}-1$ (counted with multiplicity according to Milnor
number) theoretically possible.

If other types of singularities are considered, more interesting
things happen. Invariants change, and also the type of the surfaces.
There exist  classifications of quintics with
isolated singularities (at least for those of general type), 
and interesting such examples occur for
$\tE_8$  singularities (normal form:  $z^2=y^3+\la x^2y^2+x^6$).
We will however restrict ourselves only to ordinary triple points
(locally $f(x,y,z)=0$ with $f$ a smooth plane cubic). It is
easy to see (cf.~Section \ref{section:invariants}), 
that a triple point decreases $c_1^2$
by three and $c_2$ by nine (and thus $\chi$ by one). Furthermore
the adjoint system consists of planes passing through the triple
points. We can thus establish the following table, where $\nu$
denotes the number of triple points.
$$
\halign{\qquad\quad$#$\quad&\quad\hfil$#$\hfil\quad&\quad$#$\quad&\quad$#$%
\quad&\quad$#$\quad&\quad$#$\quad&\quad#\quad\hfil\cr
\nu&c_1^2&c_2&\chi&p_q&q&Type of surface\cr
\noalign{\vskip 0.5em}
0&5&55&5&4&0&general type\cr
1&2&46&4&3&0&general type\cr
2&-1&37&3&2&0&elliptic blown up once\cr
3&-4&28&2&1&0&$K3$ blown up four times\cr
4&-7&19&1&0&0&rational\cr
5&-10&20&0&0&1&ruled over elliptic curve\cr}
$$
The most interesting thing about this table is the geometry of
the special examples, which nicely illustrates many aspects of
the general theory of surfaces. We should also note that the
constructions of surfaces are elementary, as we can impose the
triple points generically, and then simply solve  linear
equations in the coefficients. 
We  observe that the line joining any
two triple points lies by Bezout on the surface, and furthermore
is exceptional. The latter is true for any conic passing through
three triple points {\it and} lying on the quintic.

Let us now comment upon our small menagerie of surfaces, found,
so to speak in our back-yard.
\begin{description}
\item[$\nu=1\!:$] This is a double octic, the double cover of $\PP^2$
effected by projection from the triple point. (In fact all the
quintics with triple points, can be considered as double octics,
with one less triple point.) Note that not all double octics
are of this type.
\item[$\nu=2\!:$] An elliptic surface blown up. The elliptic fibration
is given by the planes through the line joining the two triple
points. The intersection consists of the line and a quartic with
two double points at its intersection with the line. The resolution
of such plane quartics (a resolution effected by the desingularisation
of the triple points) is clearly elliptic. The canonical divisor
will coincide with the elliptic fibration (plus the exceptional
divisor). The elliptic curves which arise from desingularisation,
will be bi-sections of the elliptic fibration.
\item[$\nu=3\!:$]  A $K3$ surface blown up four times. The canonical
divisor is given by the plane through the three triple points.
The intersection of the quintic with that plane is given by the
three lines of the corresponding triangle, and the residual conic,
passing through them all. If you make the resolved surface
minimal, you get a $K3$ surface with three elliptic curves $E_i$
all passing through a common point $p$ and any two $E_i$, $E_j$
($1\leq i\neq j\leq3$) also intersecting in another point $p_k$
($i$, $j$, $k$ distinct integers). Each of those elliptic curves will
give rise to elliptic fibrations, which are amusing to identify.
All the $K3$ surfaces will have Picard number at least three,
and the configuration of the elliptic curves will give a common
sublattice to them all. Conversely starting with a $K3$ surface with such a
sublattice of its Picard group, we choose a point $p$ and elliptic curves
 $E_i$ in each of the corresponding pencils, passing through $p$. 
Those will define intersection points $p_k=E_i\cap E_j$ 
($i$, $j$, $k$ distinct).
Blowing up the points to exceptional divisors $F_k$ and $G$ (the latter 
corresponding to the blow up of $p$) we can in fact write down the divisor $H$
of degree 5, effecting the birational embedding. Namely
$$
H=E_1+E_2+E_3-F_1-F_2-F_3-2G
\;.
$$
It is straightforward to check that $H^2=5$, $H\cdot E_i=0$ and 
$H\cdot F_i=1$ while $H\cdot G=2$.
Independent examples of such $K3$ surfaces are furnished by quartics with three
lines, each of which gives rise to an elliptic fibration. It could be a mildly
challenging exercise for the reader to see how such a quartic can be 
transformed into a quintic with three triple points.
It is also amusing to count parameters: there
are 56 quintic monomials, 10 conditions for a triple point at
a fixed location, and 6-dimensional family of projective linear
transformations, fixing three points. This makes 19. On the other hand $K3$
surfaces with Picard number at least three, make up a 17-dimensional family,
and we add 2 for the position of the point.

\item[$\nu=4\!:$] A rational surface blown up many times (15 or
16 times). You expect an infinite number of exceptional curves
of the first kind.
Ten of those are obvious, given by the six edges of the tetrahedron
spanned by the triple points, and the four residual conics corresponding
to each face. It is a challenge to find others.
\item[$\nu=5\!:$] A ruled surface over an elliptic curve, blown
up ten times. The twenty exceptional divisors, which come in
pairs, are obvious. Each pair of triple points determine an exceptional
line, and dually the three remaining triple points a residual
exceptional conic. This also gives a clue to the ruling. Through
six points one may always find a unique twisted cubic. For each
point on the surface, consider the twisted cubic through it and
the five triple points. By Bezout this curve has to lie on the
surface. The resolved elliptic curves will be sections. By blowing
down the ten lines (or the ten conics) we get a minimal ruled
surface. This will turn out to be the one coming from the stable
rank-two bundle on an elliptic curve.
\end{description}

Note that six or more triple points is an impossibility. Such
a putative surface would necessarily have $\chi<0$ and hence
be ruled over a curve of genus $g>1$. In such surfaces there
is no space for elliptic curves, as they can neither surject
on to the base, nor squeeze into the fibres.

Quintics with at least four triple points can be simplified
using a birational transformation, known
as {\em reciprocal transformation}
\cite[VIII \S\ 4]{semple}.
The ordinary plane Cremona transformation using the linear
system of conics through three points can be described in
suitable coordinates by the formula
$(x \cn y\cn z) \mapsto (1/x \cn 1/y\cn 1/z)$.
This formula generalises
to higher dimensions.
In particular, the space transformation
\begin{equation*}
    (x\cn y\cn z \cn w)
    \mapsto \left(\frac{1}{x}\cn\frac{1}{y}\cn\frac{1}{z}
    \cn\frac{1}{w}\right)  
\end{equation*}
simultaneously blows up the vertices and
blows down the faces of the coordinate tetrahedron.
The vertices are called {\em fundamental points} of
the reciprocal transformation.
Let $Y\subset\PP^3$ by a surface of degree $d$
not containing any of the coordinate planes. Let $m_1$, \dots, $m_4$ be the
multiplicities of $Y$ in the fundamental points.
Then the image $Y'$ of $Y$ is a surface of degree $3d-m_1-\cdots-m_4$.
In many cases $Y'$ will be singular in the fundamental points
with singularities obtained from contracting the
intersection curves of $Y$ with the coordinate planes.

For a quintic $X$ the following happens:
taking four triple points as fundamental 
points  we transform the surface into one of degree $3\cdot 5
- 4\cdot 3 = 3$. Conversely, given four points on a cubic 
surface we obtain a surface of degree $3\cdot 3 - 4 =5$ with
four new singularities, which are ordinary triple points if
the tetrahedron (spanned by the four points) cuts out smooth curves.   
This argument shows that five triple points are
maximal, and realisable by starting with a cubic cone.
The construction also allows you to find the blown up rational quintic with
four triple points. As we all know the cubic can be thought of $\PP^2$ blown
up six times. The four vertices of the tetrahedron will provide four more
blow ups, and finally the six edges of the tetrahedron, each intersect the
cubic in a residual point, each of which is blown up. This makes a total
of sixteen.

It is amusing to continue. Setting $H$ to be the hyperplane section of the 
cubic (thus $H^2=3$) and letting $E_i$ denote the exceptional curves associated
to the four vertices, and $E_{ij}$ the residual intersections associated to 
the six edges, we can write down the linear system on the cubic, which gives 
the quintic, as 
$$
3H-2\sum_iE_i-\sum_{ij}E_{ij}\;.
$$
This linear system blows down the four elliptic curves 
$$
H-E_i-E_j-E_k-E_{ij}-E_{ik}-E_{jk}
$$
each of which has self-intersection $3-6=-3$. Furthermore the six exceptional 
curves $E_{ij}$ are mapped onto lines, namely the edges of the tetrahedron 
spanned by the triple points, while the exceptional curves $E_i$ are mapped
onto residual conics. Each of the 27 lines of the cubic, 
which does not pass through a vertex,  maps onto a twisted cubic 
passing simply through all four triple points.
As there are exceptional curves of arbitrary high degree on the cubic blow-up,
we see that on a quintic there may be exceptional curves of arbitrary high
degree.

The above analysis is very elementary, and parts of it has not
too surprisingly already appeared in the literature \cite{gallarati}. It can
be generalised in a number of different ways. One, which we have
already mentioned, is to consider more general singularities
(in this way interesting surfaces of general type can be constructed 
\cite{yang}), 
the other is to consider triple points on surfaces of higher degree.

\section{Invariants}
\label{section:invariants}
A surface singularity
$P\in X$ is an {\em ordinary triple point of $X$\/}
if there exist local coordinates  
$x$, $y$ and $z$ centred at $P$
such that $X$ is given by the equation
\begin{equation*}
   x^3+y^3+z^3+\lambda xyz=0
\end{equation*}
for a $\lambda\in\CC$ with $\lambda^3\neq -27$. Such a triple point
is also called a singularity of type $\tE_6$ (resp.~$P_8$, $T_{3,3,3}$).
The minimal resolution of an ordinary triple point is given as follows:
let $\pi_P\colon\tX\rightarrow X$ be the blowup of $X$ in $P$ and
let $E_P$ be the exceptional divisor. This is a smooth elliptic
curve, given in suitable homogeneous coordinates
$(x\cn y\cn z)$ by the equation
$x^3+y^3+z^3+\lambda xyz=0$. 
In particular $\tX$ is smooth in
every point of $E_P$ and the self intersection
of $E_P$ on $\tX$ is $-3$.

Now let $X\subset \PP^3$ be a projective surface of degree $d$
with $\nu$ isolated triple points. 
Let $\sS=\left\{P_1,\ldots,P_\nu\right\}=X_{\it sing}$ be the
singular locus of $X$ and let $\tX$ be the blow-up of $X$ in
$\sS$; it is  a smooth model
of $X$, however not minimal in general. We will
denote the minimal model of $X$ by $\oX$. 
Moreover let $E=\sum_{i=1}^\nu E_i$ be the sum of all
exceptional divisors.

There are basically two types of invariants of $\tX$.
To start with, there are invariants of local nature which take
into account the number, but not the position of the
triple points. This are the Chern numbers $\conetX$, $\ctwotX$ and
the holomorphic Euler characteristic $\chitX$.
Second, there are invariants which are also influenced
by the position of the triple points. Amongst them are
the geometrical genus $\pgtX$,
the irregularity $\qtX$,
the Betti numbers $b_i\left(\smash{\tX}\right)$,
the Hodge numbers $h^{p,q}\left(\smash{\tX}\right)$ and
the Kodaira dimension $\kodtX$.
We are going to compare these invariants with the invariants
of a smooth hypersurface $X_s$ of degree $d$.

\begin{description}
\boldmath
\item[The canonical class and $\conetX$:] \unboldmath we have 
  $\KtX\sim_{lin}\piast\KXs-E$.
  This follows from the adjunction formula:
  we know that $E_i\cdot (\KtX+E_i)=0$, so if 
  $\KtX\sim_{lin}\piast K_X-\alpha E$, then $E_i\cdot (\KtX+E_i)=
  E_i\cdot (\piast\KXs-\alpha E +E_i)= 3(\alpha-1)$.
  As ${K_X}^2={\KXs}^2$ we find
\begin{equation}
  \conetX = {\KtX}^2 = {\KXs}^2-3\nu.
\end{equation}
Every triple point
diminishes $\conetX$ by three.

\boldmath\item[The Euler number $\etX=\ctwotX$:]\unboldmath 
the drop in the Euler number is
purely local and can be computed by means of topological
considerations from the Milnor number
$\mu\left(\smash{\tE_6}\right)=8$ and the Euler numbers
$e\left(E_i\right)=e\left(\mbox{\it torus\/}\right)=0$ and
$e\left(\mbox{\it point\/}\right)=1$ as follows:
\begin{align*}
  \etX &= \eXs + \nu\left(-e\left(\mbox{\it point\/}\right)
    +e\left(E_i\right)
    -\mu\left(\smash{\tE_6}\right)\right)\\
  &= \eXs - 9\nu.
\end{align*}
So every triple point diminishes the Euler number by nine.

\boldmath\item[The holomorphic Euler characteristic $\chitX$:]\unboldmath 
the Noether formula on
$\tX$ says $\chitX=\left(\conetX+\ctwotX\right)/12$, hence
\begin{equation*}
  \chitX = \chiXs -\nu.
\end{equation*}
Every triple point diminishes the holomorphic Euler characteristic by one.
\end{description}

A quick (cheating) way of computing the above `drop' in invariants 
exploits the local character of the Chern invariants. 
So consider a smooth cubic, whose invariants are $c_1^2=3$ and $c_2=9$. 
A cubic with a triple-point is a cone over an elliptic curve, 
the resolution is a ruled surface over an elliptic curve, and 
hence $c_1^2=c_2=0$ giving indeed the `drops' three and nine, deduced above.

Now we come to the other invariants.
The adjoint linear system on $X$ is cut out
by those surfaces of degree $d-4$
which pass through every point of $\sS$. Intuitively speaking,
every triple point in general position puts a linear condition on
the adjoint linear system of $X$, so it will diminish the
geometric genus $\pgtX$ by one (if not already zero).
Let $\alpha\in\NN_0$ be the discrepancy defined by
\begin{equation*}
    \pgtX = \pgXs -\nu +\alpha.
\end{equation*}
As a consequence
\begin{equation*}
    \qtX = \qXs + \alpha.
\end{equation*}
Both $\tX$ and $X_s$ are K\"ahler, so we have the Hodge decompositions
of $H^i(\tX,\ZZ)\otimes\CC\cong H^i(\tX,\CC)$ and
$H^i(X_s,\ZZ)\otimes\CC\cong H^i(X_s,\CC)$. From
the equalities $h^{p,q}=h^{q,p}=h^{4-p,q}$ and
$\bzerotX = \bzeroXs = 1$ one easily computes
\begin{align*}
    \btwotX  &= \btwoXs -9\nu +4\alpha, \\
    \honeonetX &= \honeoneXs -7\nu +2\alpha.
\end{align*}
The other Betti numbers and Hodge numbers do not give more
information:
$\bthreetX=\bonetX=2\qtX$,
$h^{1,0}\brtX=h^{0,1}\brtX=h^{2,1}\brtX=h^{1,2}\brtX=\qtX$,
$h^{2,0}\brtX=h^{0,2}\brtX=\pgtX$
and $h^{0,0}\brtX=h^{2,2}\brtX=1$.

We list the invariants in terms of $d$, $\nu$ and $\alpha$
in the following table.
\begin{table}[H]
    \centering
    \begin{math}
  \renewcommand{\arraystretch}{1.2}
    \begin{array}{|c||c|c|}\hline
    & X_s & \tX \\\hline\hline
    c_1^2   & d(d-4)^2 & d(d-4)^2-3\nu \\\hline
    c_2     & d(d^2-4d+6) & d(d^2-4d+6)-9\nu \\\hline
    \chi    & d(d^2-6d+11)/6 & d(d^2-6d+11)/6-\nu \\\hline
    p_g     & \binom{d-1}{3} & 
              \binom{d-1}{3}-\nu+\alpha\\\hline
    q       & 0 & \alpha \\\hline
    b_2     & d^3-4d^2+6d-2 & d^3-4d^2+6d-2 -9\nu +4\alpha \\\hline
    h^{1,1} & d(2d^2-6d+7)/3 & d(2d^2-6d+7)/3 -7\nu +2\alpha \\\hline
    \end{array}
    \end{math}
    \caption{The invariants of $X_s$ and $\tX$}
\end{table}
The Kodaira dimension $\kodtX$ measures the growth
of the plurigenera $P_n(\tX)=h^0(\tX,\cO_{\tX}(n\KtX))$ as $n$ grows.
In general we have  \cite[ch.~I, thm.~7.2]{barth}
\begin{equation*}
    \kodtX\left\{\begin{array}{c@{\quad}l}
        \geq 0 & \text{if $\pgtX\geq 1$ and}\\
        \geq 1 & \text{if $\pgtX\geq 2$.}
        %
    \end{array}\right.
\end{equation*}
For surfaces of general type
$P_2(\tX)={\KtX}^2+\chitX+\epsilon$, 
where $\epsilon$ is the number of exceptional divisors.
Generally we have
$P_n(\tX)=\frac12n(n-1) ({\KtX}^2+\epsilon)+\chitX$.
\section{Bounds for the number of triple points}
The surface $X$ can have only a finite number of ordinary triple points,
the maximal number depending on its degree $d$. Let $\mtd$ be the maximal
number of ordinary triple points of a degree $d$ surface.
We immediately find
\begin{equation*}
    \mu_3(1)=\mu_3(2)=0,\quad 
    \mu_3(3)=\mu_3(4)=1 \quad
    \text{and}\quad
    \mu_3(5)=5.
\end{equation*}
The only cubic surface with an ordinary
triple point is the cone over a plane smooth elliptic curve.
A quartic surface with two triple points is necessarily
singular along the line joining these two points.
As we have seen in the preceding section,
for quintics the result is due to
Gallarati \cite{gallarati}.

The position of the triple points cannot be too special,
as also the maximal number of triple points of $X$ on a given curve
or surface is bounded.
Let $C\subset\PP^3$ be a curve of degree $c$ and $V\subset\PP^3$
a surface of degree $v$.
\begin{lemma}\label{lemma:position}\noindent
    \begin{itemize}
        \item[\rm 1)] $C$ contains at most $c(d-1)/2$ triple points
                  of $X$ \lh(with multiplicity\rh).
        \item[\rm 2)] If $V$ and $X$ do not have a common component, then
                  $V$ contains at most $vd(d-1)/6$ triple points
                  of $X$ \lh(with multiplicity\rh).
    \end{itemize}
\end{lemma}
\proof Consider the linear system
$\sL_p$ of polar surfaces of $X$, i.e.~the linear system
generated by the partial derivatives of the degree $d$
polynomial defining $X$. Then $\sS$ is exactly the base locus
of $\sL_p$ and the general member $X_p\in\sL_p$ is a degree $d-1$
surface which is smooth except ordinary double points
in the triple points of $X$.
So $X_p$ does not contain a component of $C$.
In every triple point $P\in\sS$ the intersection
multiplicity of $C$ and $X_p$ in $P$ is
$\mult_P(C,X_p)\geq 2$ and thus $2\nu\leq C \cdot 
X_p = c(d-1)$. This proves 1).

The surfaces $V$, $X$ and $X_p$ intersect in a finite number of points.
In every triple point $P\in\sS$ the intersection
multiplicity of $V$, $X$ and $X_p$ in $P$ is
$\mult_P(V,X,X_p)\geq 6$. Hence $6\nu\leq V \cdot X\cdot
X_p = vd(d-1)$ and 2) holds.
\proofend

We will now discuss three bounds for $\mtd$ with 
$d\geq 6$: the polar bound, the Miyaoka bound
and the spectrum bound.

\medskip
{\bf The polar bound:} 
Suppose that $\pgtX\geq 1$. Taking a general
adjoint surface $V=K_p$ we find using  lemma
\ref{lemma:position} 2) 
\begin{equation*}
    \nu\leq \frac{1}{6}\,d(d-1)(d-4).
\end{equation*}
The condition $\pgtX\geq 1$ is satisfied for $d>6$.
This can be seen as follows.
Substitute
$\alpha=\pgtX-\binom{d-1}{3}+\nu$, then the
inequalities $1+\nu\leq\honeonetX$ and
$\qtX\geq 0$ imply 
\begin{equation*}
    \pgtX\geq \frac{1}{24}\left(d-1\right)\left(2d^2-16d+21\right).
\end{equation*}
So for $d\geq 7$ we have even $\pgtX\geq 2$.
The bound $\nu\leq \frac{1}{6}d(d-1)(d-4)$, which we
call polar bound,  even holds for $d=6$
in case $\pgtX=0$. Then
$\btwotX=\honeonetX$ and the equation
$1+\nu\leq\honeonetX$ gives
\begin{equation*}
    \nu\leq\left\lfloor\frac{1}{18}(d-1)(d^2+d-3)\right\rfloor
    = 10.
\end{equation*}
We get the following table.
%
\begin{table}[H]
    \centering
    \begin{math}
    \begin{array}{c|cccccccc}
    d       & 5  & 6  & 7  & 8 &  9 & 10 &  11 &  12 \\\hline
    \nu\leq{} & 6 & 10 & 21 & 37 & 60 & 90 & 128 & 176 
    \end{array}
    \end{math}
    \caption{The polar bound}
\end{table}

\medskip
{\bf The Miyaoka bound:} Miyaoka's famous bound \cite{miyaoka} applies only to
quotient singularities of surfaces with nonnegative Kodaira
dimension. However
there is a generalisation by Wall \cite[cor.~2]{wall} which also applies to
log-canonical singularities on surfaces with nonnegative Kodaira dimension.
Applied to triple points, we get the bound
\begin{equation*}
    \nu \leq \frac{2}{27}d(d-1)^2.
\end{equation*}
As $\pgtX\geq 2$ for $d=\deg(\tX)\geq 7$ this bound holds for
every surface with triple points of degree $\geq 7$ and we get the following
table.
\begin{table}[H]
    \centering
    \begin{math}
    \begin{array}{c|cccccccc}
    d       &  7 &  8 &  9 & 10 & 11 &  12 \\\hline
    \nu\leq{} & 18 & 29 & 42 & 60 & 81 & 107 
    \end{array}
    \end{math}
    \caption{The Miyaoka bound}
\end{table}

\medskip
{\bf The spectrum bound} \cite[sect.~14.3.2]{arnold}
uses the semicontinuity of the spectrum of a singularity.
Let $f\colon(\CC^{n+1},0)\rightarrow(\CC,0)$
be the germ of an isolated hypersurface singularity with Milnor number $\mu$.
Then the characteristic polynomial of the monodromy has $\mu$ eigenvalues
which are roots of unity, and the Mixed Hodge Structure 
on the cohomology of the Milnor fibre gives a way
to take logarithms. The precise definitions are not important for us now.
The spectrum is easy to compute for a function of the form
$f= x_0^{a_0}+ \ldots + x_n^{a_n}$: 
then the spectrum is the set of rational numbers (with multiplicity)
of the form
$i_0/a_0+...+i_n/a_n$ with the $i_j$ running from 1 to $a_j-1$.
Specifically we can take $a_i=d$ for all $i$,  and as the spectrum is
invariant under $\mu$-constant deformations, we have it now for any
homogeneous isolated singularity.
A projective hypersurface with isolated singularities has a
smooth hyperplane section, and the affine complement of the section
is a small deformation of the affine cone over the hyperplane
section, so a homogeneous isolated singularity.
The important property of the spectrum is its semicontinuity, 
in the sense that for every open
interval of length 1 the number of spectral numbers in it of the
singularity in the special fibre is at least the sum of the
spectral numbers in the same interval of all singularities in the
general fibre of a 1-parameter deformation (of negative degree).
We consider the spectrum as a divisor on $\QQ$. With this notation
the spectrum for $d=5$ is
\newcommand{\spn}[1]{\!\left({\textstyle\frac{#1}{5}}\right)}
\begin{equation*}
     1\spn{3}+
     3\spn{4}+
     6\spn{5}+
    10\spn{6}+
    12\spn{7}+
    12\spn{8}+
    10\spn{9}+
     6\spn{10}+
     3\spn{11}+
     1\spn{12}
\end{equation*}
The spectrum of an ordinary double point is 
$\left({\textstyle\frac{3}{2}}\right)$. As the open interval 
$\left(\frac{3}{5},\frac{8}{5}\right)$ contains
$31$ spectral numbers, a quintic surface can contain
at most 31 nodes. 

The spectrum for an $\tE_6$ is 
\begin{equation*}
    1\spn{3}+
    3\spn{4}+
    3\spn{5}+
    1\spn{6}.
\end{equation*}
The open interval $\left(\frac{4}{5},\frac{9}{5}\right)$
contains 40 spectrum numbers
of the quintic, and seven of $\tE_6$, so $\lfloor\frac{40}{7}\rfloor
= 5 $ is the spectrum bound. Analogous computations for higher
degree give the following table.
\begin{table}[H]
    \centering
    \begin{math}
    \begin{array}{c|cccccccc}
    d       & 5 &  6 &  7 &  8 &  9 & 10 & 11 &  12 \\\hline
    \nu\leq{} & 5 & 11 & 17 & 29 & 45 & 60 & 84 & 114 
    \end{array}
    \end{math}
    \caption{The spectrum bound}
\end{table}
\medskip
Putting all bounds together we arrive at the 
\begin{proposition}\label{proposition:bound}
    Let $X\subset\PP^3$ be a surface of degree $d\geq 3$ with
    $\nu$ ordinary triple points as its only singularities.
    Then $\nu$ is bounded as given by the following table.
    \begin{equation*}
        \begin{array}{c|cccccccccc}
        d       & 3 & 4 & 5 &  6 &  7 &  8 &  9 & 10 & 11 &  12 \\ \hline
        \nu\leq{} & 1 & 1 & 5 & 10 & 17 & 29 & 42 & 60 & 81 & 107 
        \end{array}
    \end{equation*}
\end{proposition}
%
%
%
We want to classify surfaces in $\PP^3$ with only
ordinary triple points.
In contrast to surfaces with only
ordinary double points, the class of
a surface can change if more triple points come into play.
The more triple points, the less nef $\KtX$ will become
and this makes the Kodaira dimension eventually drop.
The cases $d\leq 4$ being obvious, we can state the
\begin{proposition}\label{proposition:multiplicity}
    If $d\geq 7$, then $\tX$ is minimal.
    If $d=6$, then the smooth rational
    $(-1)$-curves on $\tX$ come from
    rational curves $C$ of degree $c\geq2$
    on $X$ through $2c+1$ 
    triple points \lh(with multiplicity\rh), while for 
    $d=5$ they come from curves of degree 
    $c\geq1$ on $X$ through $c+1$ 
    triple points.
\end{proposition}
\proof Let $C\subset X$ be a rational curve of degree $c$ such that
$\tC\subset\tX$ is a smooth rational $(-1)$-curve. Then
$-2=\deg\KtC={(\KtX+\tC)}|_{\tC}=
\KtX\cdot\tC+\tC^2 = c(d-4)-\mult\left(C,\sS\right)-1$.
Hence
\begin{equation*}
    \mult\left(C,\sS\right) = c(d-4)+1.
\end{equation*}
Applying lemma \ref{lemma:position} 1) we find that
\begin{equation*}
     c(d-4)+1= \mult\left(C,\sS\right) \leq \frac{1}{2}\,c(d-1)
\end{equation*}
and consequently $c(d-7)\leq -2$. This implies $d\leq 6$ and
$c\geq 2$ for $d=6$.
\proofend
\begin{corollary}\label{corollary:minimal}
    If $d\geq 7$, then $X$ is a minimal surface of general
    type.
\end{corollary}
\proof For $d\geq 7$ $\tX$ is minimal by 
proposition \ref{proposition:multiplicity}. $\KtX$ is effective,
 so by the Enriques-Kodaira classification we just have to show
$\conetX >0$. But $\conetX = d(d-4)^2-3\nu$, so
$\conetX\leq 0$ iff
\begin{equation*}
    \nu\geq\frac{1}{3}d(d-4)^2.
\end{equation*}
Playing this inequality against the Miyaoka bound gives
a contradiction.
\proofend
\section{Sextics}
\subsection{Exceptional curves}
We intend to study the $(-1)$-curves on $\tX$. 
The amazing thing is that there are severe restrictions.
The consequences of lemma \ref{lemma:position}
and proposition \ref{proposition:multiplicity}
are as follows. 
\begin{itemize}
    \item At most two triple points lie on a line $L\subset\PP^3$.
    \item At most five triple points lie in a plane $H\subset\PP^3$
          and if so, then $H\cdot X=3C$ for a smooth conic $C$, 
          and $\tC\subset\tX$ is
          a smooth rational $(-1)$-curve.

\end{itemize}
The conic through five triple points gives the simplest example
of a $(-1)$-curve
(this really does occur, 
see section \ref{section:less} for explicit examples). 
Such a conic will be called a $(-1)$-conic; by abuse of notation
we also call the curve $C\subset X$ a $(-1)$-curve, if
$\tC\subset\tX$ is an exceptional curve of the first kind.
The next possible candidates would be a twisted cubic curve through seven
triple points and a rational quartic curve through nine triple points.
Surprisingly the twisted cubic is impossible.
\begin{proposition}\label{proposition:cubic}
    At most six triple points lie on a cubic curve
    $D\subset X$.
\end{proposition}
\proof By lemma \ref{lemma:position} 1), a cubic 
curve contains at most seven triple
points. So assume that $\mult(D,\sS)=7$.
Then $D$ cannot be a plane curve and $D$ cannot split in three lines.
So either $D=C+L$ for a nondegenerate conic $C$ and a line $L$
with $\mult(C,\sS)=5$ and $\mult(L,\sS)=2$ or $D$ is a twisted
cubic curve.

In the first case let $H$ be the plane containing $C$. Then
$H\cdot X=3C$, so $L$ and $C$ meet in one point $P\in X\setminus\sS$.
But then $L\cdot X\geq 7$, so $L\subset X$. Hence
$L\subset T_{P,X}=H$, contradiction.

In the second case let $N$ be the net of quadrics with the twisted
cubic $D$ as its base locus. So the general (smooth) quadric
$Q\in N$ intersects $X$ in $S$ and a residual curve $D_Q$ of type
$(4,5)$ with double points in $D\cap\sS$.
But on $Q$ we have $D\cdot D_Q=(2,1).(4,5)=14$. Hence $D\cap D_Q=D\cap\sS$
for the general $Q\in N$.
But for every $P\in D\setminus\sS$ there exists a pencil of
quadrics $N_P\subset N$ having contact to $X$ at $P$.
This implies that $P\in D_Q$ for all $Q\in N_P$.
Now two things can happen. Either $N_P\subset N$ moves if we
move $P$ on $D\setminus\sS$ or $N_P$ is constant.
If $N_P$ moves it will sweep out a Zariski open subset of $N\simeq\PP^2$.
Then $D\cap D_Q=D\cap\sS$ cannot hold for the general
element of $N$. So $N_P$ is constant and for every $Q\in N_P$, $X$
has contact to $Q$ along $D$. But then for two
elements $Q\neq Q'\in N_P$ one has $Q\cdot Q'=2D+D'$ 
for some curve $D'$, contradiction.
\proofend
\begin{corollary}\label{corollary:cubic}
    There is no $(-1)$-curve of degree three on $X$.
\end{corollary}
In the case $\pgtX\geq1$ we find further restrictions for the number
of $(-1)$-curves and their degrees.
First we show that every $(-1)$-curve is irreducible.
Whenever the canonical divisor $\KtX$ is effective, 
any exceptional divisor $E$ is automatically a component, as 
$\KtX\cdot E=-1<0$. Therefore $E$ comes from a rational curve on $X$ of 
degree $c$ (by Proposition \ref{proposition:multiplicity} 
through $2c+1$ triple
points), which is contained in the base locus of the system of quadrics
through the triple points. 
This is the adjoint system; we will call every quadric in it  a {\em canonical
surface}. 

\begin{proposition}\label{proposition:disjoint}
    Let $C$ be an irreducible $(-1)$-curve on $X$, let $K$ be a canonical  
    divisor \lh(of degree 12\rh) and $C'$ the residual curve of $C$ in $K$.
    Then the strict transform $\tC\subset \tX$ of $C$ is disjoint from
    the strict transform $\tC'$ of $C'$. 
\end{proposition}
\proof 
First suppose that $C$ is a conic. Then the residual  curve $C'$ has degree
$10$, and no component of it lies in the plane through $C$. Therefore
the intersection multiplicity $C\cdot C'$ is at most $10$. As $C'$
has multiplicity $2$ in each of the five triple points on $C$, the
intersection multiplicity  is exactly  10 and $\tC$ is disjoint from
$\tC'$.

If $\deg C \geq 4$, then 
$C$ lies on an irreducible quadric $Q$. We shall show that $\tC$ is disjoint
from $\tC'$ on the blow up of $Q$ in the points $P\in\sS$.
We first suppose that $Q$ is smooth. Then $C$ is a curve of type $(a,b)$
with arithmetic genus $p_a=(a-1)(b-1)$, and $C'$ has type $(6-a, 6-b)$,
so $C\cdot C'=6(a+b)-2ab$. Suppose that $C$ has multiplicity 3 in $\tau$
points  $P\in\sS$, multiplicity $2$ in $\delta$ points, and passes simply
though $\sigma$ points. Then $3\tau +2\delta +\sigma= 2(a+b)+1$.
As $C$ is rational we have that $3\tau +\delta \leq p_a$. 
This gives $\delta+\sigma\geq 3(a+b)-ab$.
The multiplicity
of $C\cup C'$ is three in each point  $P\in\sS$,
so $C\cdot C'\geq 2\delta + 2\sigma$.
Therefore we find 
$$
\delta+\sigma\geq 3(a+b)-ab \geq \delta+\sigma\;.
$$
So $C$ intersects $C'$ only in points $P\in\sS$ and the blow up of these 
points separates both curves.

The case that $Q$ is a quadric cone with vertex outside $\sS$ is handled
in the same way. As $C$ is smooth outside $\sS$ and there does not intersect
$C'$ we conclude that $C$ does not pass through the vertex.

Finally we investigate the case that the vertex of $Q$ is a point 
$P\in \sS$. Then $K=C\cup C'$ has multiplicity $6$ in $P$. 
Let $\ol Q$ be the blow up of $Q$ in the point $P$. 
Its Picard group is generated by $E$ and $f$, with $E^2=-1$, $E\cdot f=1$
and $f^2=0$; we have $K_{\ol Q}\sim -2E-4f$. 
The strict transform of $K$ is a curve of type $3E+12f$. Let $C$ have 
multiplicity  $m$ in $P$, then its strict transform $\ol C$ is a curve of type
$aE+(2a+m)f$, with $p_a(\ol C)=(a-1)(a+m-1)$. We have 
$\ol C'\sim (3-a)E+(12-2a-m)f$, so $\ol C\cdot \ol C'=12a+3m-2a(a+m)$.
Let $\ol C$ have $\tau$ triple, $\delta$ double and $\sigma$ simple
points in $\sS\setminus\{P\}$. Then 
$3\tau +2\delta +\sigma= 4a+m+1$, $3\tau +\delta \leq p_a(\ol C)$
and $\ol C\cdot \ol C'\geq 2\delta + 2\sigma$.
Therefore
$$
\delta+\sigma\geq 6a+2m -a(a+m)\geq \delta+\sigma+\frac 12 m\;.
$$
We conclude that $m=0$, so $C$ does not pass through $P$, and that
$\tC$ is disjoint from $\tC'$. \proofend

\begin{corollary}\label{proposition:curves}
    If $\pgtX =1$, then the degree of
    every $(-1)$-curve is one of $\{2,4,5,6,7,8\}$. 
    Moreover there are at most
    $6$  such disjoint curves. 
    If $\pgtX =2$ there are at most two $(-1)$-curves
    of degree 2 or 4 and if $\pgtX \geq 3$, there is
    at most one $(-1)$-curve of degree $2$. 
\end{corollary}
\proof 
In first case all $(-1)$-curves are contained
in the unique canonical curve of degree $12$. Moreover $\conetX=-3$
and $\pgtX\neq 0$, so $\tX$ has at least three $(-1)$-curves,
which are disjoint by proposition \ref{proposition:disjoint}.
In the second case the base locus of  the adjoint system
is a curve of degree $\leq 4$.
In the last case the base locus of of the adjoint system
is a curve of degree $\leq 3$.
Now the proposition follows because there are no $(-1)$-curves
of degree 1 or  3.
\proofend
We can now determine
$\pgtX$ for $\nu=9, 10$.
\begin{corollary}\label{corollary:nine}
    If $\nu=9$, then $\pgtX=1$.
\end{corollary}
\proof Let $\nu=9$, then $\pgtX\geq 1$. If
$\pgtX\geq 2$ corollary \ref{proposition:curves}
implies that $\tX$ has at most
two $(-1)$-curves, which contradicts $\conetX=-3$ and $\pgtX\neq 0$.
\proofend
\begin{corollary}\label{corollary:ten}
    If $\nu=10$, then $\pgtX=0$.
\end{corollary}
\proof Let $\nu=10$, then  $\conetX=-6$. So if $\pgtX>0$, then
$\tX$ contains at least six $(-1)$-curves. Then $\pgtX=1$ by
corollary \ref{proposition:curves} . The only possibility for
six $(-1)$-curves is six conics $C_1$, \dots, $C_6$ which make up
$\KtX$. Blowing down the six conics gives a minimal surface
$\oX$ with $\coneoX=0$, $\ctwooX=12$ and $\KoX=\OoX$.
But there is no such surface in the Enriques-Kodaira classification.
\proofend
In fact, an $X$ with 10 triple points and $\pgtX=1$
would have an equation of the form $h_1\cdots h_6 +q^3$ where the $h_i$
define planes. Six planes intersect in 20 triple points. 
It is possible to choose 10 of them under the
condition that no three  lie on a line, but those points never lie on a 
quadric. We found our first example of a sextic with 9 triple points
by taking $h_1\cdots h_6 +q^3$ with $q$ defining
a quadric through 9 of the 10 points, chosen as required.
\par
Now let $C\subset\tX$ be a rational quartic curve such that $\tC\subset\tX$
is a smooth rational $(-1)$-curve.
$C$ is contained in a smooth quadric surface and is either of
type $(2,2)$ or $(1,3)$.
\begin{itemize}
    \item If $C$ is of type $(2,2)$, then $C$ has one double
          point in a triple point of $X$. Moreover $C$ passes simply
          through seven other triple points. $C$ is the base locus
          of a pencil of quadrics whose general member is smooth.
    \item If $C$ is of type $(1,3)$, then $C$ is smooth and passes
          simply through nine triple points of $X$.
\end{itemize}
The case of a quartic $(-1)$-curve of type $(1,3)$ 
turns out to be impossible.
\begin{lemma}\label{lemma:quartic}
    Every quartic $(-1)$-curve on $X$ is of type $(2,2)$.
\end{lemma}
\proof Assume that $C_1$ is a rational quartic $(-1)$-curve on $X$
of type $(1,3)$. Then $C_1$ is smooth and passes simply
through nine triple points of $X$, so $\nu\in\{9,10\}$.
Moreover $C_1$ is contained in a unique smooth quadric
$Q=\{q=0\}$. We have two cases.

\boldmath
{\bf Case $\nu=9$:}\unboldmath\enspace 
Then $Q$ is the unique canonical surface
(corollary \ref{corollary:nine}).
Every $(-1)$-curve is contained
in the degree twelve curve $K=Q\cdot X$.
No five triple points lie on a plane, so there exist no
$(-1)$-conics. But $\conetX=-3$ and $\pgtX=1$, so there are
at least three $(-1)$-curves. By corollary \ref{proposition:curves},
the only possibility is three $(-1)$-curves:
the curve $C_1$ of type $(1,3)$ and two other quartic
$(-1)$ curves $C_2$ and $C_3$
of types $(3,1)$ and $(2,2)$. Both $C_1$ and $C_2$
have multiplicity one in all points of $\sS$, whereas $C_3$ misses
one triple point. This contradicts $K=C_1+C_2+C_3$ and
$\mult(K,P)=3$ for all $P\in\sS$.

\boldmath
{\bf Case $\nu=10$:}\unboldmath\enspace 
By corollary
\ref{corollary:ten} we have $\pgtX=0$,
so $Q$ passes exactly through nine triple points. If $X=\{f=0\}$,
then a general element of the pencil defined by $\alpha f+\beta q^3=0$
is a sextic with $\nu=9$ ordinary triple points and $C_1$ as $(-1)$-curve.
Hence we are done using the first case.\proofend
As a further consequence we get the useful
\begin{corollary}\label{corollary:conics}
    If $C_1$, $C_2\subset X$ are two different $(-1)$-conics,
    then $C_1$ and $C_2$ meet in two distinct
    triple points of $X$.
\end{corollary}
\proof Let $H_i$ be the plane containing $C_i$, $i=1,2$.
Then $H_i\cdot X=3C_i$, so $C_1\cap C_2\neq\emptyset$. Moreover
$H_i$ is a tangent plane to $X$ at every point of $C_i$.
Thus every point of $C_1\cap C_2$ is a singular point,
i.e.~a triple point. If there is just one such point,
then $\nu\geq 9$ and $\pgtX=11-\nu$, which contradicts
one of the corollaries \ref{corollary:nine} and
\ref{corollary:ten}.
\proofend
Assume that $X$ has $\nu=9$ triple points $P_1$, \dots, $P_9$. Let
$Q$ be the unique (irreducible) canonical quadric surface
and let $K=Q\cdot X$ be the adjoint curve.
The resolution $\tX$ has at least three disjoint $(-1)$-curves,
which all are  components of $K$.
There are two main
possibilities: either $K$ is the union of all $(-1)$-curves or not.
In the first case blowing them down   gives a minimal
surface $\oX$ with $\KoX=\OoX$,  so $\coneoX=0$ and there are exactly
three $(-1)$-curves with degrees $c_1$, $c_2$ and $c_3$. 
It follows from corollary
\ref{proposition:curves} that up to permutation
\begin{equation*}
    (c_1,c_2,c_3)\in\{(2,2,8),(2,4,6),(2,5,5),(4,4,4)\}.
\end{equation*}
In the second case we end up with an effective canonical divisor
after blowing down.
In this case up to permutation the possible degrees are
$$
\displaylines{
    (2,2,2),(2,2,2,2),(2,2,2,2,2),  (2,2,4),(2,2,2,4), 
    (2,2,5),(2,2,2,5),\cr(2,2,6),(2,2,2,6),(2,2,7),
                      (2,4,4),(2,4,5).}
$$
First we will rule out some cases. For the remaining cases,
we will give explicit examples, when making
 a tour from zero to ten triple points, studying
all possible cases.
\begin{proposition}\label{proposition:makenot_up}
    If the $(-1)$-curves do not make up $K$, then there are exactly
    three with degrees $(2,2,2)$ or  $(2,2,4)$.
\end{proposition}
\proof
We exclude all other possibilities case by case.

Suppose first that there are three $(-1)$-conics. As there are only nine 
triple points the only possibility is that the three planes containing
the conics have only one point in common and that each intersection
line contains two triple points, while the remaining three points each lie 
in only one plane. Assume now that there is a fourth $(-1)$-conic.
Its five triple points have to lie on the intersection triangle
with the first three planes with at least two triple points
in the vertices. This implies that there are three triple points
on a line, thus excluding the cases $(2,2,2,2)$ and $(2,2,2,2,2)$.
Also $(2,2,2,4)$ is not possible: a quartic $(-1)$-curve $C$ passes through
$8$ triple points, and there is at least one plane $H_i$ 
containing five of them, but $C\cdot H_i=4$. 

Case $(2,2,5)$: As $\mult(C_3,\sS)=11=9+2$ the curve $C_3$ has either one
triple point or two double points in $\sS$. But every irreducible
degree five curve on $Q$ has arithmetic genus 0 or 2.
Every double point drops the genus by at least one,
every triple point by at least three (as $C_3$ does not pass through the
vertex of $Q$ if $Q$ is singular).
The only possibility is that $C_3$ has
two double points. As $C_3$ is not a plane curve, it meets
the plane of a $(-1)$-conic simply in the five triple points.
But the two $(-1)$-conics contain together  eight triple points,
so $C_3$ can have at most one double point, contradiction.
This excludes also the case $(2,2,2,5)$.

The cases $(2,2,6)$, $(2,2,2,6)$ and $(2,2,7)$ are similar.

Case $(2,4,4)$: Both $C_2$ and $C_3$ are of type $(2,2)$ and have
a singular point in $\sS$ outside $H_1$. Counting intersection
points of $C_2$ and $C_3$ gives a count $\geq 9$ (with
multiplicity). This contradicts $C_2\cdot C_3=8$.

The case $(2,4,5)$ is similar.
\proofend

We see that under the conditions of the
proposition the degree of a $(-1)$-curve is always even.
We shall show that this also holds if the degrees sum up 
to $12$ by excluding  the case $(2,5,5)$. It 
is possible to construct a  reducible
curve on $Q$ consisting of a curve of type $(1,1)$, $(2,3)$
and $(3,2)$ with the required intersection behaviour.
So we need  a different type of argument.

\begin{proposition}\label{proposition:make_up}
    If the $(-1)$-curves make up $K$, then
    \begin{equation*}
        (c_1,c_2,c_3)\in\{(2,2,8),(2,4,6),(4,4,4)\}.
    \end{equation*}
\end{proposition}
\proof Let $(c_1,c_2,c_3)=(2,5,5)$, heading for a contradiction.
Let $P_1$, \dots, $P_5$ be the five triple points on the conic $C_1$.
On $\tX$ we have $3C_1\sim_{lin}H-E_1-\ldots-E_5$, so 
\begin{equation*}
    3(C_2+C_3)\sim_{lin} 5H-2(E_1+\cdots+E_5)-3(E_6+\cdots+E_9).
\end{equation*}
Therefore there exists a degree five surface $Y=\{g=0\}$
with multiplicity two at the  points $P_1$, \dots, $P_5$ and triple points
$P_6$, \dots, $P_9$ such
that $f=hg-q^3$. Here $Q=\{q=0\}$ is the canonical quadric and
$H=\{h=0\}$ is the plane containing $C_1=\{q=h=0\}$.
As $Y$ intersects $Q$ in two irreducible curves of odd
degree, $Y$ is itself irreducible. 
A plane through three triple points intersects $Y$ in the triangle
formed by the triple points and a residual conic through them.
If one of the five double points lies 
in the plane, the conic degenerates and has multiplicity two at
the double point. But this would imply that three triple points
of the  original sextic lie on a line.
Consider the reciprocal transformation centred in the
triple points of $Y$.
The image of $Y$ will be an irreducible 
cubic surface $Y'$ with five  points of multiplicity two. 
Therefore $Y'$ has nonisolated singularities: it has a double line.
So $Y$ itself has a double curve of degree at most three, passing through
the five points $P_1$, \dots, $P_5$. This means that the  conic $C_1$ 
is a component of the double line, so the sextic surface $X$ is singular
along $C_1$.
\proofend 
\subsection{Sextics with seven or less triple points}
\label{section:less}
For all surfaces $X$ with up to four triple points $\tX$
is minimal by proposition \ref{proposition:multiplicity}.
No three triple points lie on a line,
thus $\pgtX=10-\nu$. There are no constraints on the position of the
triple points except that no three are on a line, so we get
an equation of $X$ just by solving linear equations in the
coefficients.
If four triple points $P_1$, $P_2$, $P_3$, $P_4$ 
lie in a plane $H\subset\PP^3$, then $X\cdot H$ is a degree six
curve with four triple points and hence splits into
three conics: $X\cdot H=C_1+C_2+C_3$.
Then $X$ has an equation of the form
\begin{equation*}
    hg+q_1q_2q_3=0
\end{equation*}
with $H=\{h=0\}$ and $C_i=\{h=q_i=0\}$, $i=1,2,3$. Here $g$ is a
degree five form which vanishes in $P_1$, \dots, $P_4$ to the second order.
In any case, the base locus of the system of adjoint surfaces
consists only of the triple points.

Imposing a fifth triple point $P_5$ opens the possibility of a
$(-1)$-curve. This happens if and only if $P_1$, \dots, $P_5$ lie
on a plane $H$. Then $X\cdot H=3C$ for a nondegenerate conic $C$
and $\tC\subset\tX$ is the $(-1)$-curve. Every such surface
has an equation of the form
\begin{equation*}
    hg+q^3=0
\end{equation*}
with $H=\{h=0\}$ and $C=\{h=q=0\}$. Here $g$ is a degree five form
vanishing doubly in $P_1$, \dots, $P_5$. Then the base
locus of the system of adjoint surfaces is exactly $C$.
If the five triple points do not lie in a plane we find
$\{P_1$, \dots, $P_5\}$ as base locus. In any case $\pgtX=5$.

Let $P_6$ be another triple point. Since six triple points
cannot lie on a conic the geometric genus will drop
by one: $\pgtX=4$. We end up with the same cases
as for $\nu=5$ yielding as base locus
$C\cup\{P_6\}$ resp. $\{P_1,\ldots,P_6\}$.

Examples of sextics with $\nu\leq 6$ triple points and
base locus $\{P_1,\ldots,P_{\nu}\}$ of the adjoint system
can be given as follows. Let $Q_i=\{q_i=0\}$ be generators
of the linear system of quadrics through $\{P_1,\ldots,P_{\nu}\}$,
$i=0$, \dots, $10-\nu$. Then the general element of the linear
system spanned by the mixed third powers of the
$q_i$ has only triple
points in $P_1$, \dots, $P_\nu$.
For $\nu\leq5$ every surface is of this form. For $\nu=6$ the linear 
system of mixed third powers has dimension $\binom 63-1=19$ while the system
of all sextics with triple points has dimension $\binom 96-1-6\cdot10=23$.

Now we go for a seventh triple point $P_7$. Again the geometric genus
drops: $\pgtX=3$. The base locus cannot be a degree three curve
by proposition \ref{proposition:cubic}. If $P_1$, \dots, $P_5$ lie on a conic
$C$, the base locus is $C\cup\{P_6,P_7\}$. 
If not, the
net of quadrics defined by $P_1$, \dots, $P_7$
has a zero-dimensional base  of the form
$\{P_1,\ldots,P_7,P\}$ for a eighth point $P\in\PP^3$, which may be
infinitely near to one of the points $P_1$, \dots, $P_7$.
Now the mixed third powers of the $q_i$ 
have an additional singularity in $P$. They
form a system of dimension $9$, which is
four less than the dimension of the system of all sextics with triple
points. To find an equation of such a surface it suffices to give
one possibly reducible sextic not passing through $P$; the surface obtained
by adding a general combination of third powers has then only $7$
triple points. As such an extra sextic we can take the product $g_1g_2$ of
a cubic $g_1$ with nodes in
$P_1$, \dots, $P_4$  passing simply through $P_5$, $P_6$, $P_7$
and a cubic $g_2$ with nodes in $P_5$, $P_6$, $P_7$
passing simply through $P_1$, \dots, $P_4$.

In all cases considered so far we end up with
$\conetX>0$ and $\pgtX\geq 3$, so
$\tX$ is a surface of general type.
\subsection{Eight triple points}
We distinguish the sextics with eight triple points by their
geometric genus and their $(-1)$-curves. 

We can always choose seven of the eight points so that no five
lie on a conic. These seven triple points
$P_1$, \dots, $P_7$
(no three on a line, no five on a conic) determine a net
of quadrics spanned by $Q_i=\{q_i=0\}$, $i=1,2,3$.

\medskip
{\bf The case $\pgtX=3$:} this means that 
the eighth triple point $P_8$ is the eighth base point of
the net.  For a general ternary cubic form $f_3$ the surface
$X=\{f_3(q_1,q_2,q_3)=0\}$ is a sextic with only triple
points in $P_1$, \dots, $P_8$.
Here $\pgtX=3$, thus $\qtX=1$  and $\tX$ is minimal because $K$ is 
effective and has no fixed components,
so $\tX$ is minimal properly elliptic. 
This elliptic surface is fibred over an elliptic curve, 
namely the plane elliptic curve given by $f_3=0$. 
The fibration is induced by the net, i.e. given by $(q_1,q_2,q_3)$ 
and each fibre is the base locus of a pencil of quadrics, 
in a way the reader easily can work out. Note that the elliptic resolutions
of the triple points are sections of this fibration, and hence all isomorphic 
to the base. (This can also be seen 
by noting that the linear parts of the $q_i$'s at the 
basepoints are linearly independent.)

\medskip
{\bf The case $\pgtX=2$:}
we assume that $P_8$ is not a base point of the net. Let the
pencil of quadrics through $P_1$, \dots, $P_8$ be spanned by
$Q_1$ and $Q_2$. Let $C$ be the base locus of the pencil.
For the $(-1)$-curves we can have four
different cases:
\begin{itemize}
    \item one $(-1)$-curve of type $(2,2)$,
    \item two $(-1)$-conics,
    \item one $(-1)$-conic and
    \item no $(-1)$ curves at all.
\end{itemize}

\medskip
{\bf One $(-1)$-curve of type $(2,2)$:}
in this case the equation of $X$ has 
a very special form.
\begin{lemma}\label{lemma:equation}
    If $\tX$ contains a quartic $(-1)$-curve $C$ of type $(2,2)$, then
    the equation of $X$ has the form
    \begin{equation*}
       q_0g+q^3=0
    \end{equation*}
    with $Q_0=\{q_0=0\}$  a quadric cone through eight triple points
    with vertex in one of them  and
    $Q=\{q=0\}$  a smooth quadric through the eight triple points.
    The pencil of quadrics with base locus $C$ is spanned by
    $Q$ and $Q_0$.
    Moreover $Y=\{g=0\}$ is a quartic surface passing through the vertex
    of $Q_0$ with seven double points in the other seven triple points.
\end{lemma}
\proof Let $P_1\in\sS$ be the double point of $C$ and let
$P_2$, \dots, $P_7$ be the other triple points on $C$.
Let $M$ be the pencil of quadrics with base locus $C$.
The general element $Q\in M$ is smooth and intersects $X$ in $C$ and
a residual curve $C_Q$ of type $(4,4)$ passing simply through $P_1$ and
doubly through $P_2$, \dots, $P_7$. Since $C$ and $C_Q$ do not have
a common component we have $C\cap C_Q = \{P_1,\ldots,P_8\}$
in view of $C\cdot C_Q=(2,2)\cdot(4,4)=16$. 
Now fix a point $P\in C\setminus\sS$.
There exists a $Q_0\in M$ which has contact to $X$ at $P$. In particular
$P\in C_{Q_0}$, implying that $C$ and $C_{Q_0}$ have a common component. 
Hence $X\cdot Q_0=2C+C'$ for a curve $C'$ of type $(2,2)$.
Thus $\mult(X\cdot Q_0,P_1)\geq4$, so $Q_0$ is singular in $P_1$.
Hence $\mult(C',P_1)\geq 2$ and $X\cdot Q_0=3C'$ by Bezout.
In view of $C\subset Q_0$ the quadric $Q_0$ has to be
a quadric cone with vertex $P_1$.
Then $X$ has an equation as demanded.
\proofend

The condition that a quadric 
$\sum \la_i q_i$
in the net of quadrics through 7 points
(in general position) 
is singular is that  there exists a point $P$ in which all derivatives
of $\sum \la_i q_i$ vanish. Eliminating the $\la_i$ gives that the 
maximal minors of 
$$
\def\part#1#2{\frac{\partial#1}{\partial#2}}
\def\partx#1{\part {q_{#1}}{x} & \part{q_{#1}}{y} &
           \part {q_{#1}}{z} &\part {q_{#1}}{w}}
 \left|
\matrix \partx 1 \\ \partx 2 \\ \partx 3\endmatrix
\right |
$$
have to vanish. This locus is known as the 
Steiner curve of the net.
For $Q$ we take a smooth quadric
through the eight points. 
There exists a six parameter
family of quartics through $P_8$ with only
nodes in $P_1$, \dots, $P_7$ (quadrics in the net give a five parameter
family, but the general quartic is not of that form).
Let $Y=\{g=0\}$ be a general such quartic,
then $X=\{q_0g+q^3=0\}$ is a sextic with eight triple
points and one $(-1)$-curve $C$ of type $(2,2)$.
Altogether this construction has 13 moduli:
we find seven independent sextics with triple points in the
given points, and the configuration of points
has seven moduli.
The minimal model $\oX$ of $X$ has $\coneoX=1$, so
$\tX$ is of general type. 

Such a surface has a characterisation as the minimal model of
a double cover of $F_2$ (the resolution of a quadric cone), branched along a 
5-section disjoint from the minimal section (the node) and that section. 
If $S$ denotes a section with $S^2=2$, and $F$ a fibre, the branch curve is 
simply $5S+(S-2F)$. As $K=-2B$ we have $K+B=S-F=(S-2F)+F$ with $S-2F$ as 
fixed component. Thus the free part of the canonical pencil gives a genus 
two fibration.

A generic quadric in the pencil intersects the sextic in a $(6,6)$ curve 
which splits up into the fixed component (the $(2,2)$ curve) and a residual 
$(4,4)$ curve with seven double points, passing simply through the node 
$P_8$ of the $(2,2)$ curve. The genus of the desingularisation is 
indeed $3\times 3-7=2$, giving us the genus two pencil.

In fact if $\{q+tq_0=0\}$ is a quadric in the pencil, 
the residual curve is the intersection $\{q+tq_0=tq^2+g=0\}$.
The canonical system on such a genus two curve is given by 
adjunction as the pencil of $(2,2)$ curves passing through the seven nodes. 
Thus the involution on the surface will be given as follows. For each point 
$P$ choose $t$ such $q+tq_0$ vanish at $P$, then consider the residual 
intersection with the base locus of the pencil of quadrics through $P_1$, 
\dots, $P_7$, $P$ and the quartic $\{tq^2+g=0\}$.

\medskip
{\bf One $(-1)$-conic:} such a sextic can be obtained as
the reciprocal transform of the previous sextic in
$P_1$, $P_2$, $P_3$ and $P_8$ (which necessarily do not 
lie in one plane). The quadric $Q$ transforms into a smooth quadric
$Q'=\{q'=0\}$, $Q_0$ transforms into a
plane $H_0'=\{h_0'=0\}$ and $Y$ transforms into
a quintic $Y'=\{g'=0\}$ with three triple points and
five double points. So the transform $X'$ of $X$
satisfies the equation $h_0'g'+{q'}^3=0$. The image
of $C$ is a $(-1)$-conic lying in the plane $H_0'$.
Conversely, every such sextic can be transformed into
one with a $(-1)$-curve of type $(2,2)$: just take
as fundamental points the three triple points not on $H_0'$ and
a fourth triple point on $H_0'$. The base locus of
the pencil of quadrics consists of the $(-1)$-conic
and another conic not contained in $X'$. This family
again has 13 moduli.

\medskip
{\bf Two $(-1)$-conics:} the sextics $X$ with two $(-1)$-conics 
$C_1$ and $C_2$ are easily identified
as those satisfying an equation of the form
$h_1h_2g+q^3=0$. Here $H_i=\{h_i=0\}$ are planes containing
$C_i$, $i=1,2$. By corollary \ref{corollary:conics}, the two conics
intersect in two triple points, say $P_1$ and $P_2$.
Then $Y=\{g=0\}$ is a quartic surface through $P_1$ and
$P_2$ with only nodes in $P_3$, \dots, $P_8$ and $Q=\{q=0\}$
is a general element of the pencil. There are ten
linearly independent sextics with the eight triple
points: four from the pencil spanned by
$h_1h_2$ and $q$, another four of the form
$h_1^2h_2^2q_1$ where $q_1$ is one of the four quadrics
through $P_3$, \dots, $P_8$ and finally $h_1h_2^2k_1$,
$h_1^2h_2k_2$ for cubics $k_1$ and $k_2$ through
all eight points such that $k_i$ has double points
in the triple points not contained in $h_i$, $i=1,2$.
The point configuration has five moduli,
so we get a 14 parameter family. Again every such
sextic is of general type.

The minimal models will have $p_g=2$ and $c_1^2=2$ and come equipped with 
genus three fibrations. Those will be defined by, 
in analogy with the case of a $(-1)$ curve of type $(2,2)$, 
as residual $(4,4)$ curves, with nodes at $P_3\dots P_8$, 
passing through $P_1$, $P_2$. We leave it to the reader to work out 
the details.

\medskip
{\bf No $(-1)$-curves:}
it is not so immediate  to construct such surfaces.
We follow the classical   construction for
asyzygetic eight-nodal quartics \cite{cayley,rohn}, which we now
recall.
Seven points $P_1$, \dots, $P_7$ in general
position define a net of quadrics spanned by
smooth quadrics $Q_i=\{q_i=0\}$, $i=1,2,3$.
All quartics with nodes in $P_1$, \dots, $P_7$ are given by
$g+f_2(q_1,q_2,q_3)=0$, where $f_2$ is a ternary quadratic form and
$g$ is a fixed nodal quartic, which we can take as the product of
a cubic with nodes in
$P_1$, \dots, $P_4$ passing simply through $P_5$, $P_6$, $P_7$
and a plane through $P_5$, $P_6$, $P_7$.
A sufficient condition for an eighth singular point is that the 
derivatives $dg$, $dq_1$, $dq_2$ and $dq_3$ are linearly dependent.
The vanishing of the determinant
\begin{equation*}
    \begin{vmatrix}
        \frac{\partial g}{\partial x} &
        \frac{\partial q_1}{\partial x} &
        \frac{\partial q_2}{\partial x} &
        \frac{\partial q_3}{\partial x} \\[3pt]
        \frac{\partial g}{\partial y} &
        \frac{\partial q_1}{\partial y} &
        \frac{\partial q_2}{\partial y} &
        \frac{\partial q_3}{\partial y} \\[3pt]
        \frac{\partial g}{\partial z} &
        \frac{\partial q_1}{\partial z} &
        \frac{\partial q_2}{\partial z} &
        \frac{\partial q_3}{\partial z} \\[3pt]
        \frac{\partial g}{\partial w} &
        \frac{\partial q_1}{\partial w} &
        \frac{\partial q_2}{\partial w} &
        \frac{\partial q_3}{\partial w}
    \end{vmatrix}
\end{equation*}
defines a degree six surface $\Delta$, called the Cayley dianode surface.
It is
the closure of the
locus of points $P_8$ such that there
exists a quartic surface with only nodes in
$P_1$, \dots, $P_8$. In general the 
dianode surface has triple
points in $P_1$, \dots, $P_7$ as its only singularities, but it can reducible:
if four points lie in a plane this plane becomes a component.

We now look at sextics with seven triple points in general position.
We find 14 such sextic equations:
the ten third powers of the three quadrics
    $q_1$, $q_2$ and $q_3$,
three equations of the form $\{gq_i=0\}$ for a quartic
    $g$ which is {\em not} a conic of the three quadrics and finally
a sextic $\delta$ not of this form, for which one can take the
Cayley dianode surface.
The vanishing of all second derivatives in an eighth point
$P_8$ gives ten linear equations in the 14 coefficients.
Note however that the columns of the coefficient matrix are not 
linearly independent.
We observe that the all second derivatives of a third
power have the function itself as factor:
$$
(f^3)_{ij}=3f(2f_if_j+ff_{ij})\;.
$$
For a product we have
$$
(f^2g)_{ij}= (2f_if_j+ff_{ij})g+f(2f_ig_j+2f_jg_i+f_{ij}g+fg_{ij})
$$
and a similar expression for the derivatives of $fgh$.
Suppose that $f\neq0$. Then we can divide all partials $(f^3)_{ij}$
by $f$ and 
after subtracting the same multiple of $2f_if_j+ff_{ij}$ 
from all partials $(f^2g)_{ij}$ we can again divide by $f$. 
{}From $(fg^2)_{ij}$ we can get $2g_ig_j+gg_{ij}$.
So if one of the quadrics does not
vanish  (i.e., the point is not the 8th base
point of the net), the 10 columns of the third powers give at most 6
independent ones and our matrix reduces to a
$10\times 10$ matrix.
One obtains a determinant of degree 28, but the rank
of the matrix is eight on the dianode surface, whose equation is 
double factor of the determinant. We are left with a surface
$\Delta'$ of degree 16. It has as double curves the 21 lines joining the
triple points, the 7 cubics though 6 of the 7 points and the Steiner curve
of the net. The dianode surface $\Delta$ intersects $\Delta'$ exactly
in its singular locus.

Now we get two possibilities to construct a sextic, by taking
$P_8$ in $\Delta$ or in $\Delta'$.
If $P_8$ is a general point of the Cayley dianode surface, 
there exists a quartic surface $Y=\{g=0\}$ with
only nodes in $P_1$, \dots, $P_8$. Then a general linear combination
of $gq_1$, $gq_2$, $q_1^3$, $q_1^2q_2$, $q_1q_2^2$, $q_2^3$
defines a sextic $X$ with only triple points in
the eight points and containing the base locus of the
pencil spanned by $Q_1$ and $Q_2$.
The surface will be a (minimal) elliptic surface, 
whose elliptic pencil is given by
the pencil of quadrics through the points $P_1$, \dots, $P_8$:
the residual intersection on a smooth quadric
of the pencil is a curve of type $(4,4)$ with $8$
double points.
The point configuration having eight moduli,
we get a 13 parameter family.

Choosing $P_8\in\Delta'$ general we obtain a sextic with eight triple points
not containing parts of the base locus of the pencil.
The number of parameters is now $8+4=12$.
The surface is again a (minimal) elliptic surface, with elliptic pencil 
given by the pencil of quadrics: a general intersection is a
$(6,6)$-curve with eight triple points, hence the geometric genus is one.
\subsection{Nine triple points}
Assume that $X$ has $\nu=9$ triple points $P_1$, \dots, $P_9$. Let
$Q$ be the unique (irreducible) canonical quadric surface
and let $K=Q\cdot X$ be the adjoint curve.
By propositions \ref{proposition:make_up} and \ref{proposition:makenot_up}
the resolution $\tX$ has exactly three disjoint $(-1)$-curves 
$C_1$, $C_2$, $C_3$.
If $C_1+C_2+C_3=K$,
blowing down $C_1$, $C_2$ and $C_3$  gives a minimal
surface $\oX$ with $\coneoX=0$, $\ctwooX=24$ and
$\KoX=\OoX$. Then $\tX$ is a $K3$ surface blown up in
three points. 
Otherwise we end up with an effective canonical divisor
after blowing down $C_1$, $C_2$ and $C_3$. 
This implies $\kodtX=1$ and thus $\tX$ is the blowup of
a minimal properly elliptic surface in three points. As it will
have the same basic invariants as a $K3$ surface, it will be obtained from an 
elliptic such by a series of logarithmic transforms, i.e. making some elliptic
 fibres multiple.
\subsubsection*{The sextic $K3$ surface}
Using corollary \ref{corollary:conics}, the multiplicities of the
$(-1)$-curves in the nine triple points are easily found
to be (up to a permutation of triple points) as in the following table.
\begin{table}[H]
    \centering
    \begin{tabular}{|l|l||c|c|c|c|c|c|c|c|c|}\hline
       type &
       &$P_1$&$P_2$&$P_3$&$P_4$&$P_5$&$P_6$&$P_7$&$P_8$&$P_9$\\\hline\hline
       \multirow{3}{1.3cm}{$(4,4,4)$} 
        & $C_1$ & 2 & 1 & 0 & 1 & 1 & 1 & 1 & 1 & 1 \\\cline{2-11}
        & $C_2$ & 0 & 2 & 1 & 1 & 1 & 1 & 1 & 1 & 1 \\\cline{2-11}
        & $C_3$ & 1 & 0 & 2 & 1 & 1 & 1 & 1 & 1 & 1 \\\hline\hline
       \multirow{3}{1.3cm}{$(2,4,6)$} 
        & $C_1$ & 1 & 0 & 0 & 0 & 0 & 1 & 1 & 1 & 1 \\\cline{2-11}
        & $C_2$ & 0 & 2 & 1 & 1 & 1 & 1 & 1 & 1 & 1 \\\cline{2-11}
        & $C_3$ & 2 & 1 & 2 & 2 & 2 & 1 & 1 & 1 & 1 \\\hline\hline
       \multirow{3}{1.3cm}{$(2,2,8)$} 
        & $C_1$ & 1 & 0 & 0 & 0 & 0 & 1 & 1 & 1 & 1 \\\cline{2-11}
        & $C_2$ & 0 & 1 & 1 & 1 & 0 & 0 & 0 & 1 & 1 \\\cline{2-11}
        & $C_3$ & 2 & 2 & 2 & 2 & 3 & 2 & 2 & 1 & 1 \\\hline\hline
    \end{tabular}
    \caption{multiplicities of the $(-1)$-curves}
    \label{table:multK3}
\end{table}

\begin{proposition}
    If $(c_1,c_2,c_3)=(4,4,4)$, then $X$ satisfies
    an equation of the form
    \begin{equation*}
        q_1q_2q_3+q^3=0\;,
    \end{equation*}
    where $Q=\{q=0\}$ is the unique canonical surface. Each $Q_i=\{q_i=0\}$
    is a quadric cone through eight triple points with vertex
    in one of them and $C_i=Q\cdot Q_i$, $i=1$, $2$, $3$.
\end{proposition}
\proof The three quartic $(-1)$-curves
are of type $(2,2)$ by lemma \ref{lemma:quartic}.
Lemma \ref{lemma:equation} guarantees the existence
of three quadric cones $Q_i=\{q_i=0\}$ such that
$X\cdot Q_i=3C_i$, $i=1,2,3$. Thus
$X\cdot(Q_1+Q_2+Q_3)=3(C_1+C_2+C_3)$, hence a equation
of $X$ and $q_1q_2q_3$ differ by the cube of a quadratic
polynomial vanishing in all triple points.
Since $\pgtX=1$, such
polynomial defines the unique canonical surface.
\proofend

Let $Q_i=\{q_i=0\}$ be the quadric cones
with vertices $P_i$ such that
$Q_i$ passes through $P_{i+1\bmod 3}$ but not through
$P_{i+2 \bmod 3}$, $i=1,2,3$.
We choose $P_1$, $P_2$ and  $P_3$ as 
$(1\cn0\cn0\cn0)$, $(0\cn1\cn0\cn0)$ and $(0\cn0\cn1\cn0)$.
If we require that the quadrics also pass through 
$(0\cn0\cn0\cn1)$ and $(1\cn1\cn1\cn1)$ their equations
can be written (inhomogeneously) as
\begin{align*}
    q_1 &= z^2+a_1y+b_1z-(a_1+b_1+1)yz,\\
    q_2 &= x^2+a_2z+b_2x-(a_2+b_2+1)xz,\\
    q_3 &= y^2+a_3x+b_3y-(a_3+b_3+1)yx,
\end{align*}
where the $a_i$ and $b_i$ are constants.
The quadrics  intersect in a zero-dimensional scheme 
of length $8$ containing at least six distinct points. 
To find the canonical surface $Q=\{q=0\}$ we have to pick
six points $P_4$, \dots $P_{9}$.
It is easier to specify the remaining
scheme of length two. For a Zariski open set of the stratum
it will consist of two points in general position, which we take as
$(0\cn0\cn0\cn1)$ and $(1\cn1\cn1\cn1)$. 
Thus we require that $Q$ does not pass through these two points.
To compute its equation $q=0$ we note that
the $q_i$ lie in the ideal defining $(1\cn1\cn1\cn1)$; they can be written
in vector form as 
\begin{equation*}
    \begin{pmatrix}
        q_1 \\
        q_2  \\
        q_3  
    \end{pmatrix} =
     \begin{pmatrix}
        0          & (b_1+z)z  &  (a_1-z)y \\
        (a_2-x)z   & 0         &  (b_2+x)x \\
        (b_3+y)y   & (a_3-y)x  &  0 
    \end{pmatrix}
    \begin{pmatrix}
        1-x \\
        1-y  \\
        1-z  
    \end{pmatrix}
    \;.
\end{equation*}
In the triple points of $X$ the $q_i$ vanish, while $(x,y,z)\neq(1,1,1)$.
So the determinant of the matrix vanishes.
Dividing it by $xyz$ gives the
inhomogeneous equation 
\begin{equation*}
    q=(a_1-z)(a_2-x)(a_3-y)+(b_1+z)(b_2+x)(b_3+y)
\end{equation*}
which is indeed the sought after quadric (note that the two terms
$xyz$ cancel).
The general element of the pencil
$\alpha q_1q_2q_3+\beta q^3$ defines a sextic $X$ with
nine triple points with multiplicities as 
in the first part of table \ref{table:multK3}.

A particular example is obtained by taking $b_i=0$,
$a_i=-1$, so $(q_1,q_2,q_3)=(z^2-y,x^2-z,y^2-x)$
and $P_{3+i}=(\eta^{4i},\eta^{2i},\eta^i)$ with $\eta$ a primitive
seventh root of unity. The number of parameters in the
construction is 7 (the $a_i$, $b_i$  and $(\alpha\cn\beta)$). 

Now consider the reciprocal transformation with fundamental points
$P_1$, $P_2$, $P_4$ and $P_5$. 
In general the image $X'$ of $X$ will have nine ordinary
triple points as only singularities. It can be checked that this is
indeed the case for the particular example, where we take 
$(\eta^{4},\eta^{2},\eta)$ and $(\ol\eta^{4},\ol\eta^{2},\ol\eta)$
as third and fourth fundamental point. 
Then the image $H_1'$ of $Q_1$ is a plane,
the image $Q_2'$ of $Q_2$ is a quadric cone and
the image $Y_3'$ of $Q_3$ is a cubic surface with four nodes.
The image $Q'$ of $Q$ being again
a quadric, the sextic $X'$ is given by an equation of the form
\begin{equation*}
    h_1'q_2'g_3'+{q'}^3=0.
\end{equation*}
Note that $C_1$, $C_2$ and $C_3$ are mapped onto 
$(-1)$-curves $C_1'$, $C_2'$ and $C_3'$ of degree $2$, $4$ and $6$.

Surfaces with an equation of this type form a seven dimensional family.
The nodes of the cubic, the vertex of the quadric and two points
in the plane  can be specified, while the last two are then to be found
among the four other intersection points of the cubic, the
cone and the plane. The dimension 
of the Zariski tangent space to the
equisingular stratum  in the specific example is  $15+7$, so we obtain a
full seven parameter family of sextics with nine triple
points and $(c_1,c_2,c_3)=(2,4,6)$.

Now let $X$ have $(c_1,c_2,c_3)=(2,4,6)$ and an
equation of the form $h_1q_2g+q^3=0$, where $h_1$ determines a
plane $H_1$, $q_2$ a quadric cone $Q_2$ and $g_3$ a four-nodal cubic $Y_3$.
Consider
the reciprocal transform with fundamental points
$P_2$, $P_5$, $P_6$ and $P_7$ as in table
\ref{table:multK3}. 
Continuing with the specific example above, one can take the
images of the points $(\eta,\eta^4,\eta^2)$ and 
$(\ol\eta,\ol\eta^4,\ol\eta^2)$ as $P_6$ and $P_7$. The reciprocal image
has nine ordinary triple points, so again this property holds on an
open dense set in the parameter space of all sextics
with $(c_1,c_2,c_3)=(2,4,6)$.
Here the image $H_1'$ of $H_1$ is again a plane, the
image $H_2'$ of $Q_2$ is also a plane and the image
$Y'_4$ of $Y_3$ is a quartic surface with one
triple point and six double points.
The image $Q'$ of $Q$ is again a quadric,
so the image $X'$ of $X$ is given by an equation of the form
\begin{equation*}
    h_1'h_2'g'_4+{q'}^3=0.
\end{equation*}
This time the curves $C_1$, $C_2$ and $C_3$ are mapped onto 
$(-1)$-curves $C_1'$, $C_2'$ and $C_3'$ of degrees
$2$, $2$ and $8$.

Once more we check that this is 
a full seven parameter family of sextic surfaces with nine triple
points and $(c_1,c_2,c_3)=(2,2,8)$.
A direct construction of the family starts with seven points in general
position, of which we  choose
one as the triple point for the quartic $Y_4$ and divide the 
remaining six into two groups of three, each determining a plane
$H_i$. The intersection $Y_4\cap H_1 \cap H_2$ consists of
four points. Two of them can be triple points for a sextic
(remember that three are not allowed on a line) in the 
pencil $\alpha q^3 + \beta h_1h_2g_4$. 
An explicit example starts 
coordinate vertices and the points 
$(\la:1:1:1)$, $(1:\la:1:1)$ and $(1:1:\la:1)$.
We get 
$$
\displaylines{
h_1= t,\qquad h_2=x+y+z-(\la+2)t \cr
\qquad g_4=\la(\la+1)(\la+2)(xy+xz+yz)^2-(2\la+1)^2(\la+1)xyz(x+y+z)
 \hfill \cr \hfill
 -\la(2\la+1)(xy+xz+yz)(x+y+z)t+(2\la+1)^2(\la+2)xyzt\;.\qquad\cr} 
$$
The intersection line $H_1\cap H_2$ is now a double
tangent of $Y_4\cap H_1$ so we take the  two points of tangency
as last two triple points. We find
$$
q=(2+\la)(xy+xz+yz)-(2\la+1)(\la+1)(x+y+z)t \;.
$$

There is a certain amount of choice in the fundamental points
of the reciprocal transformations. They depend also on the position
of the triple points. So we can conclude that we obtain correspondences
between our families of sextics with triple points, whose exact nature 
we did not determine. We will say that our families are
`related via reciprocal transformations'.

We summarise our findings.
\begin{theorem}\label{theorem:K3}
    For every $(c_1,c_2,c_3)\in\{(2,2,8),(2,4,6),(4,4,4)\}$ there
    exists a seven parameter family of sextic surfaces with
    nine triple points and three $(-1)$-curves of degrees
    $c_1$, $c_2$ and $c_3$. The three families are
    related via reciprocal transformations.
    For every such surface $X$ its minimal desingularisation
    $\tX$ is a $K3$ surface blown up in three points.
    Moreover $X$ satisfies an equation of the form
    \begin{equation*}
        \begin{array}{l@{\quad}l}
            q_1q_2q_3+q^3=0 & \text{if $(c_1,c_2,c_3)=(4,4,4)$.}\\[1ex]
            h_1q_2g_3+q^3=0 & \text{if $(c_1,c_2,c_3)=(2,4,6)$,}\\[1ex]
            h_1h_2g_4+q^3=0 & \text{if $(c_1,c_2,c_3)=(2,2,8)$.}
        \end{array}
    \end{equation*}
    Here $Q=\{q=0\}$ is the unique canonical surface.
    The three exceptional curves are in each case obtained
    as the locus where one of the three forms in the product
    and $q$ vanish.
    In the last case $g_4$ defines a
    quartic surface with a triple point
    six double points.
    In the second case, $g_3$ defines a four nodal cubic.
\end{theorem}
It may now be amusing to digress on the geometry of the $K3$ surfaces obtained.
Consider the case $(4,4,4)$. Let $E_i$ be  the image in the minimal $K3$ 
surface of the exceptional curve in the
resolution of $P_i$. Notice that 
$E_i^2=2$ if $i\leq3$ and $E_i^2=0$ otherwise. Furthermore for $i\neq j$ 
we have $E_i\cdot E_j=2$ if $i,j\leq 3$ while $E_i\cdot E_j=3$ otherwise.
For future reference let us denote the case $E^2=0$ as the first type 
and $E^2=2$ as the second type, and by slight abuse of terminology, 
also the corresponding  triple points.

It is now easy to write down 
$$
2H=\sum_iE_i-8(C_1+C_2+C_3)
$$
and from the above it is straightforward to check that $(2H)^2=24$.
Notice also that $\sum_i E_i$ is an even divisor.

Conversely given such a configuration of curves $E_i$ in the Picard
group we need to choose the points $c_i$ to be blown up carefully.
Any divisor $E$ with $E^2=2$ on a $K3$ surface determines a net, and thus
a Jacobian curve of the net, corresponding to the curve-locus of singular 
points of singular members. Each $E_i$, $i\leq3$ determines such a
curve $J_i$. The point $c_1$ has to be chosen on $J_1$. The point $c_2$ 
has to lie on the intersection of $J_2$ with some element in $|E_1|$ 
singular at $c_1$ and we expect only a finite number of such choices.
Furthermore $c_3$ has to lie on the intersection of $J_3$ with some 
element of $|E_2|$ singular at $c_2$, and finally some element of $|E_3|$
singular at $c_3$ should pass through $c_1$. This indicates that $c_1$ should 
be chosen with care, 
and you thus expect only a finite number of configurations 
of points $c_i$. Further conditions are that the remaining $E_i$
also pass through the points $c_i$.
We expect a 11-dimensional family of $K3$ surfaces 
with the appropriate sublattices, and therefore four independent conditions.
Their exact nature remains mysterious.

To consider the two remaining cases, one writes down respectively
$$
2H=\sum_iE_i-4C_1-8C_2-12C_3
$$
and
$$
2H=\sum_iE_i-4C_1-4C_2-16C_3
$$
with different intersection matrices, and similar conditions on the
points $c_i$; in the last case there will also be a third type  of
elliptic curve  ($E_5^2=6$).

To see how those are related to the first case we note (once again) the fact 
that the intersection of a sextic with a plane passing through exactly
three triple point gives an elliptic curve $F$ with $F^2=-3$ after the
resolution. Given a tetrahedron of triple points, thus means replacing 
the exceptional vertices with the elliptic curves of the faces.

To make this explicit return to the first case $(4,4,4)$. Let us denote 
by $i$, $j$, $k$ different integers strictly greater than three, and $m$ an 
integer less or equal to three. We get four cases for $F$ namely
$$
F=H+3(C_1+C_2+C_3)-E_i-E_j-E_k
$$
or
$$
F=H+3(C_1+C_2+C_3)-E_1-E_2-E_3
$$
with $F\cdot C_m=1$ and also
$$
F=H+4C_1+2C_2+3C_3-E_1-E_j-E_k
$$
or
$$
F=H+4C_1+3C_2+2C_3)-E_1-E_2-E_k
$$
with $F\cdot C_1=0$, $F\cdot C_2=2$, $F\cdot C_3=1$ and 
$F\cdot C_1=0$, $F\cdot C_2=1$, $F\cdot C_3=2$ respectively.

Thus in the first two cases we get elliptic curves $F$ of the first 
type, while in the last two cases, curves of the second type.

In the case $(4,4,4)$ we have three curves of the second type and
six of the first. If we choose a tetrahedron with no triple points of 
the second type, or three, the number of elliptic curves of first or
second type will not change. However if there are one or two triple 
points of the second type, after the transformation the number of each 
type will change from $6,3$ to $4,5$ landing us in the case $(2,4,6)$.
We leave it to the reader to continue the analysis and show how we can 
get from situation $(2,4,6)$ back to $(4,4,4)$ or to $(2,2,8)$.
\subsubsection*{The sextic properly elliptic surface}
Now we turn to the remaining cases where $C_1+C_2+C_3$
does not make up $K$.
Using corollary \ref{corollary:conics} again, the multiplicities of the
$(-1)$-curves in the nine triple points are easily found
to be (up to a permutation of triple points) as in the following table.
\begin{table}[H]
    \centering
    \begin{tabular}{|l|l||c|c|c|c|c|c|c|c|c|}\hline
       type &
       &$P_1$&$P_2$&$P_3$&$P_4$&$P_5$&$P_6$&$P_7$&$P_8$&$P_9$\\\hline\hline
       \multirow{3}{1.3cm}{$(2,2,2)$} 
        & $C_1$ & 0 & 0 & 1 & 1 & 1 & 1 & 1 & 0 & 0 \\\cline{2-11} 
        & $C_2$ & 1 & 1 & 0 & 0 & 1 & 1 & 0 & 1 & 0 \\\cline{2-11} 
        & $C_3$ & 1 & 1 & 1 & 1 & 0 & 0 & 0 & 0 & 1 \\\hline\hline 
       \multirow{3}{1.3cm}{$(2,2,4)$} 
        & $C_1$ & 0 & 0 & 1 & 1 & 1 & 1 & 1 & 0 & 0 \\\cline{2-11} 
        & $C_2$ & 1 & 1 & 0 & 0 & 1 & 1 & 0 & 1 & 0 \\\cline{2-11} 
        & $C_3$ & 1 & 1 & 1 & 1 & 0 & 1 & 1 & 1 & 2 \\\hline\hline 
    \end{tabular}
    \caption{multiplicities of the $(-1)$-curves}
    \label{table:multelliptic}
\end{table}
Let us construct a sextic $X$ with nine triple points and
$(c_1,c_2,c_3)=(2,2,2)$. Then $X$ has an equation of the form
$f=h_1h_2h_3g+q^3=0$. So take three planes in general position
$H_i=\{h_i=0\}$, $i=1,2,3$ and let $L_1=H_1\cap H_2$,
$L_2=H_1\cap H_3$ and $L_3=H_2\cap H_3$ be the lines of
intersection.
On each line $L_i$ choose two different points $P_{2i-1}$, $P_{2i}$.
Finally choose
points $P_{i+7}\in H_i$ not on the lines.
The nine points
determine a unique quadric $Q=\{q=0\}$. Let $H_4=\{h_4=0\}$
be the plane through $P_7$, $P_8$ and $P_9$.
The reducible cubic $qh_4=0$ is an element of the pencil
of cubics through all points with double points
in $P_7$, $P_8$ and $P_9$.
Let $Y=\{g=0\}$ be another such cubic. Then 
the general element of the net
\begin{equation*}
    \alpha q^3 + \beta h_1h_2h_3qh_4 + \gamma h_1h_2h_3 g
\end{equation*}
has  only isolated triple points at the nine points.

Note that the canonical divisor will be a $(3,3)$ curve with three nodes,
the intersection of $Y$ with $Q$.
This is an elliptic curve and constitutes a 
reduced multiple fibre $F_0$. As $K=(m-1)F_0$ we conclude $m=2$. 
(In fact we have, as $X$ has the same invariants as a $K3$ surface, 
$K=\sum_i(m_i-1)F_{m_i}$, where all the multiple fibres $F_{m_i}$ are 
fixed components of the canonical divisor) Thus 
we have blown up a so called `fake $K3$' 
surface, obtained from an elliptic $K3$ surface, by making a fibre double. 
Thus in particular $P_2=2$ as $2K=F$.
In fact the bicanonical divisors are given by quartics with nodes at the
nine triple points. We can write down two independent such quartics, namely
$q^2$ and $h_1h_2h_3h_4$. 
The pencil spanned by those cuts out the residual to  
the configuration of the three exceptional conics (with multiplicity two).
As we have noted $q^2$ cuts out the double fibre $F_0$, 
while the other quartic 
will cut out the elliptic curve $F$ (intersection with the plane $h_4$) 
along with the three exceptional divisors. Thus $F$ is a general fibre, which
has been blown up three times.

It could be interesting to investigate the position of the resolution 
curves $E_i$ (corresponding to the triple points $P_i$). 
Three of those $E_7$, $E_8$, $E_9$ span the plane $h_4$ and 
thus each of them intersects $F$ three times. They also intersect $F_0$ 
twice (corresponding to its nodes). Furthermore each of them intersects 
exactly one of the exceptional divisors, namely the one, whose plane the
corresponding triple point happens to lie on. Thus in the minimal model 
these curves are elliptic four-sections, thus self-intersection $-2$. The
remaining six intersect $F_0$ simply, hence are bi-sections, meaning
self-intersection $-1$. 
Each of those intersects exactly two exceptional divisors 
and in the minimal model they will intersect $F$ twice rather than being 
disjoint from it, as in the resolution.

The reader could easily work out the intersection matrix of those nine curves
and from that write down the hyperplane section. A similar analysis could 
be made for the second case (to be discussed below), 
and the relation between the two via the reciprocal 
transformation, elucidated by considering the possible tetrahedrons of triple 
points.

In order to search for a tenth triple point we need
explicit equations. After a change of
coordinates we may assume that the planes are the
sides of the coordinate tetrahedron. The remaining 
coordinate transformations are given by diagonal matrices.
We take
$P_7=(0\cn1\cn\lambda\cn0)$,
$P_8=(\mu\cn0\cn1\cn0)$ and
$P_9=(1\cn\nu\cn0\cn0)$. A cubic through these points and 
also through $(0\cn0\cn0\cn1)$ is
\begin{align*}
    g =\, & w^2(a_1x+a_2y+a_3z)  \\
      & +w( a_4x(\nu x - y - \mu\nu z) 
           +a_5y(\la y - \la\nu x - z)
           +a_6z(\mu z - x - \mu\la z))\\
      & +(\nu x- y - \mu\nu z)(\la y - \la\nu x - z)(\mu z -x -\mu\la z),
\end{align*}
while a quadric through $P_7$, $P_8$ and $P_9$ is
\begin{align*}
    q =\, & b_4b_5b_6w^2
            + w(b_1b_4x+b_2b_5y+b_3b_6z) \\
        & + b_4x(\nu x - y - \mu\nu z)
          + b_5y(\la y - \la\nu x - z)
          + b_6z(\mu z - x - \mu\la z),
\end{align*}
where the parameters are chosen with hindsight.
The condition that $Y$ and $Q$ intersect the coordinate
axes in the same points leads to easy equations between the
coefficients. We get
\begin{align*}
    g =\, & w^2(\la\nu b_5b_6x
               +\la\mu b_4b_6y
               +\mu\nu b_4b_5z) \\
        &{} + w\big( \la b_1x(\nu x - y - \mu\nu z)
            +    \mu b_2y(\la y - \la\nu x - z)
            +    \nu b_3z(\mu z - x - \mu\la z)\big)\\
        &{} +(\nu x - y - \mu\nu z)
           (\la y - \la\nu x - z)
           (\mu z - x -\ mu\la z)\;.
\end{align*}
The formulas contain nine parameters, three of which can be 
removed by coordinate transformations. The family depends
therefore on eight moduli.

We can specialise to a symmetric equation by taking
$\la=\mu=\nu=1$, $b_1=b_2=b_3=b$ and $b_4=b_5=b_6=\sqrt{a}$.
With $a=b=1$ the surface $X=\{xyzg+q^3=0\}$
is a sextic with only nine
triple points.

In this example we can compute the tangent
space to the equisingular stratum to have dimension $15+8$.
This shows that  our construction fills up a whole component.
Note that the multiplicities of the $(-1)$-curves are just
as in table \ref{table:multelliptic}.

Consider the reciprocal transformation centred in
$P_6$, $P_7$, $P_8$ and $P_9$. The image $H_1'$ of
$H_1$ is a plane, the image $H_2'$ of $H_2$ is a plane,
the image $Q_3'$ of $H_3$ is a quadric cone and the image
$Q_4'$ of
$K$ is  a smooth quadric. So the image $X'$ of $X$ is given
by an equation 
\begin{equation*}
    \alpha {q'}^3 + \beta h_1'h_2'q_1'q'+ \gamma h_1'h_2'q_1'q_2' =0\;.
\end{equation*}
The images $C_1'$, $C_2'$ and $C_3'$ of
$C_1$, $C_2$ and $C_3$ are $(-1)$-curves
of degrees $2$, $2$ and $4$.
As before we see that
the two cases
$(c_1,c_2,c_3)=(2,2,2)$ and $(c_1,c_2,c_3)=(2,2,4)$ are
related by  reciprocal transformations.
So there exists an eight parameter family of sextics
with nine triple  points and $(c_1,c_2,c_3)=(2,2,4)$.

Now  the canonical divisor will be a smooth $(2,2)$ curve. 
It is again a fibre occurring with multiplicity two.

We summarise:
\begin{theorem}\label{theorem:elliptic}
    For every $(c_1,c_2,c_3)\in\{(2,2,2),(2,2,4)\}$ there
    exists an eight parameter family of sextic surfaces with
    nine triple points and three $(-1)$-curves of degrees
    $c_1$, $c_2$ and $c_3$. The two families are related
    via reciprocal transformations.
    For every such surface $X$ its
    minimal desingularisation $\tX$ is a minimal properly
    elliptic surface blown up in three points.
    Moreover $X$ satisfies an equation of the form
    \begin{equation*}
        \begin{array}{l@{\quad}l}
            h_1h_2h_3g+q^3=0 & \text{if $(c_1,c_2,c_3)=(2,2,2)$,}\\[1ex]
            h_1h_2q_3g+q^3=0 & \text{if $(c_1,c_2,c_3)=(2,2,4)$.}
        \end{array}
    \end{equation*}
    Here $Q=\{q=0\}$ is the unique canonical surface.
    The three exceptional curves are in each case obtained
    as the locus where the three forms in the product
    and $q$ vanish.
    In the first case $g$ defines a
    cubic with three double points.
    In the second case, $g$ defines a smooth quadric.
\end{theorem}
\subsection{Ten triple points}
Let $X$ be a sextic with $\nu=10$ triple points.
As for the type of $\tX$ in the classifi\-cation, we have the
\begin{proposition}
    If $\nu=10$, then $X$ is rational.
\end{proposition}
\proof We will show that $\kodtX=-\infty$, then the result follows
from the Enriques-Kodaira classification.

Assume that $\kodtX=0$. Then $\tX$ would be an Enriques
surface $\oX$ blown up in six points. Hence $\ptwotX=1$, so there
exists a unique quartic bicanonical surface $Y$ intersecting
$X$ in a degree 24 curve $D$ made up by the six
$(-1)$-curves $C_1$, \dots, $C_6$ of degrees $c_1$, \dots, $c_6$
(remember $2\KoX=\OoX$). On the one hand we get
$\mult(D,\sS)=2(c_1+\ldots+c_6)+6=54$ from proposition
\ref{proposition:multiplicity}.
On the other hand $2\KtX\sim_{lin}4H-2E$, hence
$\mult(D,\sS)=10\cdot2\cdot3=60$, contradiction.

Now assume that $\kodtX\geq 1$, then $\ptwotX\geq 2$.
So we have at least two quartic bicanonical
surfaces $Y_1$ and $Y_2$ intersecting in a degree 16 curve $D$.
We must have $\mult(D,\sS)=10\cdot 2\cdot 2=40$.
There exists a decomposition
$D=D_1+D_2$ with $D_1\subset X$ and no component of $D_2\neq 0$
is contained in $X$. Let $d_i=\deg D_1$, $i=1,2$.
Then $\mult(D,\sS)=\mult(D_1,\sS)+\mult(D_2,\sS)\leq
5d_1/2+2d_2<40$, contradiction.
\proofend
Every sextic with ten triple points
is a specialisation of a family of sextics with
nine triple points: if
$Q=\{q=0\}$ is the unique quadric through  nine out of the ten
triple points, the general element of the pencil
$\alpha f+\beta q^3=0$ is a sextic with nine triple points,
where $f$ is a defining equation for  $X$.
A sextic with $\nu=10$ is likely to be found in any of the five
families with nine triple points described above, as a triple point
gives seven conditions. However the equations on the coefficients become rather
formidable,  and we have only succeeded in one case by imposing extra
symmetry. 
We start with the first family of properly elliptic surfaces
with equations
\begin{equation*}
    \alpha q^3 + \beta xyzqw + \gamma xyz g
\end{equation*}
where the planes $H_1=\{x=0\}$, \dots, $H_4=\{w=0\}$
are the faces of the coordinate tetrahedron 
and the cubic and quadric are as given above.
We use the remaining freedom in coordinate transformations to
place the putative tenth triple point in $(1\cn1\cn1\cn1)$. 
We compute in the affine chart $w=1$.

The condition for a triple point is then that the function, its derivatives
and the second order derivatives vanish at $(1,1,1)$.
This gives ten equations which are linear in 
$\alpha$, $\beta $ and $\gamma$, so we may eliminate them: the
maximal minors of the coefficient matrix have to vanish.
We have
\begin{align*}
    \frac{\partial\, xyzg}{\partial x}     &= yzg + xyzg_x, \qquad
    \frac{\partial^2\, xyzg}{\partial x^2} = 
             2yzg_x + xyzg_{xx} \quad\text{and}\\
    & \frac{\partial^2 \,xyzg}{\partial x\,\partial y}  
       = zg +xzg_x+yzg_y+xyz g_{xy}.
\end{align*}
Now we plug in $x=y=z=1$. From $g$ we get
\begin{align*}
    &\la\nu b_5b_6+\la\mu b_4b_6+\mu\nu b_4b_5\\
    &+\la b_1(\nu - 1 - \mu\nu)
     +\mu b_2(\la - \la\nu - 1)
     +\nu b_3(\mu - 1 - \mu\la)\\
    &+(\nu - 1 - \mu\nu)
      (\la - \la\nu - 1)
      (\mu - 1 - \mu\la),
\end{align*}
an expression which we continue to denote by $g$.
We get also expressions for all derivatives.
Likewise we have 
\begin{align*}
    q =\,&   b_4b_5b_6 + (b_1b_4 + b_2b_5 + b_3b_6) \\
         & + b_4(\nu  - 1 - \mu\nu  )
           + b_5(\la  - \la\nu  - 1 )
           + b_6(\mu  - 1 - \mu\la  ) \;.
\end{align*}
Moreover
\begin{align*}
    \frac{\partial q^3}{\partial x}           & = 3q^2q_x,\qquad
    \frac{\partial^2 q^3}{\partial x^2}  = 3q^2q_{xx}+6qq_x^2 \quad\text{and}\\
    & \frac{\partial^2 q^3}{\partial x\partial y} = 3q^2q_{xy}+6qq_xq_y.
\end{align*}
All these are divisible by $q$. After dividing by $q$
our matrix has the following form:
\begin{equation*}
  \begin{pmatrix}
   q^2&  3qq_x & \ldots & 3qq_{xx}+6q_x^2 &\ldots& 3qq_{xy}+6q_xq_y & \ldots\\
   q  & q +q_x & \ldots & 2q_x+q_{xx}     &\ldots& q+q_x+q_y+q_{xy} & \ldots\\
   g  & g +g_x & \ldots & 2g_x+g_{xx}     &\ldots& g+g_x+g_y+g_{xy} & \ldots
  \end{pmatrix}
\end{equation*}
We simplify this matrix by subtracting $3q$ times the second row from the
first row to remove all second derivatives from the first row.
After that we apply only column operations.
A computation reveals that $q_{xx}+2\nu q_{xy} +\nu^2 q_{yy}=0$ and
also $g_{xx}+2\nu g_{xy} +\nu^2 g_{yy}=0$. Analogous equations
hold for the other second partials. After multiplying the column
containing $q_{xy}$ by $\nu$ and further column operations we get
\begin{equation*}
    \begin{pmatrix}
        -2q^2 & -q^2 & \ldots & 2q^2-6qq_x & \ldots & p_{xy}   & \ldots \\
          q   &  q_x & \ldots &  q_{xx}   & \ldots & 0 & \ldots \\
          g   &  g_x & \ldots &  g_{xx}   & \ldots & 0 & \ldots
    \end{pmatrix},
\end{equation*}
where 
\begin{equation*}
    p_{xy} = (\nu^2+\nu+1)q^2-3(\nu+1)q(\nu q_y+q_x)+3(\nu q_y+q_x)^2\;.
\end{equation*}
Now the entries in the columns with two zeroes have to vanish, for
otherwise $\alpha=0$ and the equation for the sextic is divisible by $xyz$.
We obtain the three equations
\begin{align*}
    (\nu^2+\nu+1)q^2-3(\nu+1)q(\nu q_y+q_x)+3(\nu q_y+q_x)^2 &= 0\,, \\
    (\mu^2+\mu+1)q^2-3(\mu+1)q(\mu q_x+q_z)+3(\mu q_x+q_z)^2 &= 0\,,\\
    (\la^2+\la+1)q^2-3(\la+1)q(\la q_z+q_y)+3(\la q_z+q_y)^2 &= 0\,.
\end{align*}
This shows that the locus we are after consists of several
components. 
The discriminant of the first equation, as a quadratic
form in $q$ and $\nu q_y+q_x$, equals $-3(\nu-1)^2$,
which implies that no solution is defined over $\RR$.
In principle $q_x$ is now expressible in terms of $\la$, $\mu$, 
$\nu$ and $q$, so the first row becomes divisible by $q^2$.
We may divide by $q^2$ because there is no quadric through all ten
triple points. 

This simplification is not enough to solve the equations.
To obtain manageable equations we impose symmetry.
We take $\la=\mu=\nu$, $b_1=b_2=b_3=b$ and $b_4=b_5=b_6=\sqrt{a}$. 
This gives 
\begin{align*}
    g =\,&  \la^2a(x+y+z) +
          b\left(x(\la x-y-\la^2z)+
                 y(\la y-\la^2x-z)+
                 z(\la z-x-\la^2y)\right) \\
    & +(\la x-y-\la^2 z)(\la y-\la^2x -z)(\la z -x-\la^2y)
\end{align*}
and after dividing by $\sqrt a$
\begin{align*}
    q =\, & a+b(x+y+z)
        +x(\la x- y - \la^2 z) + y (\la y -\la^2x - z) +z(\la z -x -\la^2y).
\end{align*}
The second derivatives evaluated in $(1,1,1)$ give  only four
different equations due to the symmetry in $x$, $y$ and $z$.
Our three equations reduce to
\begin{equation*}
    \frac14(\la-1)^2q^2+3(\la+1)^2(q_x-\frac12q)^2=0\,.
\end{equation*}
Substituting the values of $q$ and $q_x$ gives up to a constant
\begin{equation*}
    3(\la-1)^2(b+a/3-\la^2+\la-1)^2+(\la+1)^2(b+a+\la^2-\la+1)^2 =0\,.
\end{equation*}
The first seven columns of our matrix above reduce to three independent
ones. After multiplication of the last column by $(\la+1)^2$ we can
use the equation above to eliminate $q_x$. Upon division by $\la$ we get
\begin{equation*}
    \begin{pmatrix}
        -2q^2  & -q^2  & -q^2  \\
        q  & q_x & q_{xy}- q_{xx} \\
        g  & g_x & g_{xy}- g_{xx}&
    \end{pmatrix}\,.
\end{equation*}
As $q\neq 0$ we get as second  equation
\begin{equation*}
    q(g_x+g_{xx}-g_{xy})-(q_x+q_{xx}-q_{xy})g
    -2q_x(g_{xx}-g_{xy})+2(q_{xx}-q_{xy})g_x
\end{equation*}
and by subtracting a suitable multiple of the first equation it
becomes divisible by $(\la^2-1)^2$, giving as final equations
\begin{align*}
    & a^2+3ab+3b^2-3a\la = 0\,,\\
    & 3(\la-1)^2(b+a/3-\la^2+\la-1)^2+(\la+1)^2(b+a+\la^2-\la+1)^2 =0\,.
\end{align*}
These equations define a reducible curve, but no component
is defined over $\QQ$. To get a specific example note that
in characteristic $31$ there is a solution $\la=2$, $a=9$ and $b=-11$.
For these values of the parameters one finds explicit points.
A {\it Macaulay\/} computation shows that there is in fact a
unique sextic with only isolated triple points in the ten points.
This shows that for a general solution over $\CC$ of the equations above
a unique surface with ten triple points exists.

Computing the dimension of the tangent space to the
equisingular stratum for the specific example gives $15+3$.
However if one looks at the number of equations and the number of
variables it seems that in the general case the solution
space has to be one-dimensional. 

If one leaves out $P_7$, $P_8$ or $P_9$, one finds 
a pencil of sextics with nine triple points in the
remaining ones, which contains as reducible curve two conics
with five points in it and a degree eight curve
with one triple point and six nodes.
This is a surface with nine triple points and
$(c_1,c_2,c_3)=(2,2,8)$!
So the
sextic with ten triple points lies in the closure 
of several families.
We get solutions in the 
closure of other families by applying Cremona transformations.

One could also start out with the family $\alpha q^3+\beta q_1q_2q_3$.
The surfaces with ten triple points in this family are reciprocally
related with the surfaces in the other families with the same
number of parameters.
As the Cremona transformation depends on the position of the
points it might be possible that one finds a real surface in this
family.
Unfortunately the corresponding equations are too difficult to solve.

Altogether we have proved the following
\begin{theorem}\label{theorem:ten}
    For every triple
    \begin{equation*}
        (c_1,c_2,c_3)\in
        \{(2,2,2),(2,2,4),(2,2,8),(2,4,6),(4,4,4)\},
    \end{equation*}
    the closure of 
    the seven parameter family of sextics with $\nu=9$ triple points
    and $(-1)$-curves of degrees $c_1$, $c_2$ and $c_2$
    contains at least a one parameter family of rational sextics with
    ten triple points. 
\end{theorem}
\begin{corollary}
    $\mu_3(6)=10$.
\end{corollary}
\subsection{Summary}
\begin{theorem}\label{theorem:sextics}
    The sextics with triple points fall into 18 classes according
    to the following table.
    \begin{table}[H]
        \centering
        \begin{math}
        \begin{array}{|c||c|c|c|c|c|c|c|c|c|l|}\hline
\nu&c_1^2&c_2&\chi& p_g&q&b_2  &h^{1,1}&\#(-1)&\kod&\oX\\\hline\hline
0 & 24 &108 &11 & 10 & 0 & 106 & 86 & 0 &  2 & \text{general type} \\\hline
1 & 21 & 99 &10 &  9 & 0 &  97 & 79 & 0 &  2 & \text{general type} \\\hline
2 & 18 & 90 & 9 &  8 & 0 &  88 & 72 & 0 &  2 & \text{general type} \\\hline
3 & 15 & 81 & 8 &  7 & 0 &  79 & 65 & 0 &  2 & \text{general type} \\\hline
4 & 12 & 72 & 7 &  6 & 0 &  70 & 58 & 0 &  2 & \text{general type} \\\hline
5 &  9 & 63 & 6 &  5 & 0 &  61 & 51 & 1 &  2 & \text{general type} \\\hline
5 &  9 & 63 & 6 &  5 & 0 &  61 & 51 & 0 &  2 & \text{general type} \\\hline
6 &  6 & 54 & 5 &  4 & 0 &  52 & 44 & 1 &  2 & \text{general type} \\\hline
6 &  6 & 54 & 5 &  4 & 0 &  52 & 44 & 0 &  2 & \text{general type} \\\hline
7 &  3 & 45 & 4 &  3 & 0 &  43 & 37 & 1 &  2 & \text{general type} \\\hline
7 &  3 & 45 & 4 &  3 & 0 &  43 & 37 & 0 &  2 & \text{general type} \\\hline
8 &  0 & 36 & 3 &  2 & 0 &  34 & 30 & 1 &  2 & \text{general type} \\\hline
8 &  0 & 36 & 3 &  2 & 0 &  34 & 30 & 2 &  2 & \text{general type} \\\hline
8 &  0 & 36 & 3 &  2 & 0 &  34 & 30 & 0 &  1 & \text{elliptic}     \\\hline
8 &  0 & 36 & 3 &  3 & 1 &  36 & 32 & 0 &  1 & \text{elliptic}     \\\hline
9 & -3 & 27 & 2 &  1 & 0 &  25 & 23 & 3 &  1 & \text{elliptic}     \\\hline
9 & -3 & 27 & 2 &  1 & 0 &  25 & 23 & 3 &  0 & K3                  \\\hline
10& -6 & 18 & 1 &  0 & 0 &  16 & 16 &  &-\infty&\text{rational}    \\\hline
    \end{array}
        \end{math}
    \end{table}
    \noindent 
    All numbers denote  invariants of the corresponding surface and
    $\#(-1)$ denotes the number of $(-1)$-curves
    distinguishing $\tX$ from its minimal model $\oX$.
\end{theorem}
\section{Higher degree}
Surfaces of degree $d\geq7$ with many triple points are surfaces of general
type by corollary \ref{corollary:minimal}.
It is, as with surface with many ordinary double points,
very difficult to find explicit examples of high degree
with many ordinary triple points.
A septic surface ($d=7$)
can have at most 17 triple points by the spectrum bound.

We construct a one parameter family of septics
with 16 triple points. The symmetric group $S_4$ acts
on the polynomial ring $C[x,y,z,w]$ by permutation
of the variables.
The $\CC$ vector
space of all $S_4$-symmetric polynomials of degree seven has
dimension eleven and is generated by the polynomials
\begin{equation*}
    \sigma_1^7,\,
    \sigma_1^5\sigma_2,\,
    \sigma_1^4\sigma_3,\,
    \sigma_1^3\sigma_2^2,\,
    \sigma_1^3\sigma_4,\,
    \sigma_1^2\sigma_2\sigma_3,\,
    \sigma_1\sigma_2^3,\,
    \sigma_1\sigma_2\sigma_4,\,
    \sigma_1\sigma_3^2,\,
    \sigma_2^2\sigma_3,\,
    \sigma_3\sigma_4\,.
\end{equation*}
Here $\sigma_i$ denotes the $i$-th elementary symmetric
polynomial, $i=1$, \dots, $4$. 
Now take as 16 triple points the
$S_4$-orbit of length four generated by $P_1=(1\cn0\cn0\cn0)$
consisting of the vertices of the coordinate tetrahedron and an orbit of twelve
points generated by a point
$R_1=(\lambda\cn\mu\cn\nu\cn\nu)$.

For a $S_4$-symmetric septic, the condition to have a triple point in
$P_1$ implies that the coefficients of the first four polynomials vanish.
Imposing a triple point in $R_1$ gives 10 equations in 7 coefficients,
which by symmetry reduce to  seven  equations. 
This system of linear equations gives
an $7\times7$ matrix whose determinant is up to a constant
\begin{align*}
    \nu^5 & 
    (\lambda-\mu)^4
    (\lambda-\nu)^5
    (\mu-\nu)^5
    (\lambda+\nu)
    (\mu+\nu)
    (\lambda+\mu+2\nu)^4\\
    & \cdot (\lambda\mu-\nu^2)
    (2\lambda\mu+\lambda\nu+\mu\nu)
    (\lambda\mu+2\lambda\nu+2\mu\nu+\nu^2)^3.
\end{align*}
It is easily checked that all solutions except
$\lambda+\nu=0$ or equivalently $\mu+\nu=0$ correspond to
either degenerate surfaces or to orbits with less than twelve
points.
For $\lambda+\nu=0$ the orbit of triple points is generated
by $(-\nu\cn\mu\cn\nu\cn\nu)$. 
The coefficients of the symmetric polynomials are now
easily determined.
For general values of $(\mu\cn\nu)\in\PP^1$ one has indeed a septic with
16 isolated ordinary triple points. Computing the dimension of the
tangent space to the equisingular stratum gives one,
so all surfaces in the family have $S_4$-symmetry.
There is no irreducible surface with 17 triple points in this family.
\begin{theorem}
    The general element of the
    one parameter family of $S_4$-symmetric septics given by
    \begin{align*}
        (\mu-\nu)^3\nu(\sigma_1^2\sigma_2\sigma_3 &-
            \sigma_1\sigma_3^2 - \sigma_1^3\sigma_4) 
        -(\mu+\nu)\nu^3\sigma_1\sigma_2^3
           -(\mu+\nu)(\mu^3-\nu^3)\sigma_2^2\sigma_3\\
          & +(\mu+\nu)(\mu-\nu)^2(\mu+2\nu)\sigma_1\sigma_2\sigma_4
             +(\mu+\nu)(\mu-\nu)^3\sigma_3\sigma_4 =0
    \end{align*}
    has 16 ordinary triple points as its only singularities.
\end{theorem}
\begin{corollary}
    $16\leq\mu_3(7)\leq 17$.
\end{corollary}
\begin{remark}
    Whenever 16 points are invariant under the symmetric group $S_4$,
    it is tempting to ask if they form a Kummer configuration.
    This is not the case here.
\end{remark}
\subsection*{Acknowledgements}
We thank Duco van Straten for participating in our search
for sextics with many triple points.
\noindent%
\vfill

\parindent=0pt
\small
Addresses of the authors:\\
Stephan Endra\ss\ ({\tt endrass@micronas.com})\\
Micronas GmbH,
Hans-Bunte-Stra\ss e 19,
D 79108 Freiburg, Germany\\
Ulf Persson ({\tt ulfp@math.chalmers.se}), 
Jan Stevens ({\tt stevens@math.chalmers.se})\\
Matematik,
Chalmers tekniska h\"ogskola,
SE 412 96 G\"oteborg, Sweden
\end{document}